\documentclass{article}    % Enable this line and disable the
\usepackage{graphicx}          % Include this line if your
\usepackage{epsfig}
\usepackage{graphicx}
\usepackage{makeidx}
\usepackage{multicol}
\usepackage{verbatim}
\usepackage{subfigure}
\usepackage{mathrsfs}
\usepackage{snapshot}

\usepackage{amsmath}

\usepackage{crop}
\crop
\usepackage{amssymb}%,minilisp,mycltxdefs}%,amsthm}%,verbatim}
\usepackage{graphicx}

\def\E{\mathbb{E}}

\def\ba{\begin{array}}
\def\ea{\end{array}}
\def\bi{\begin{itemize}}
\def\ei{\end{itemize}}

\newtheorem{defn}{Definition}[section]

\newtheorem{thm}{Theorem}[section]
\newtheorem{lem}{Lemma}[section]

\newtheorem{cond}{Condition}[section]

\begin{document}

%----
%\begin{frontmatter}
%\runtitle{Insert a suggested running title}  % Running title for regular
                                              % papers but only if the title
                                              % is over 5 words. Running title
                                              % is not shown in output.
\begin{comment}
----
\title{Identification of Finite Dimensional Linear Systems Driven by L\'evy processes} % Title, preferably not more
                                                % than 10 words.

%\thanks[footnoteinfo]{This paper was not presented at any IFAC
%meeting. Corresponding author M.~T.~Cicero. Tel. +XXXIX-VI-mmmxxi.
%Fax +XXXIX-VI-mmmxxv.}

\author[GL]{L\'aszl\'o Gerencs\'er}\ead{gerencser.laszlo@sztaki.mta.hu},    % Add the
\author[MM]{M\'at\'e M\'anfay}\ead{manfay@sztaki.mta.hu}            % e-mail address

\address[GL]{MTA SZTAKI}  % Please supply
\address[MM]{MTA SZTAKI, Central European University}             % full addresses
%\address[Baiae]{The White House, Baiae}        % here.

\begin{keyword}                           % Five to ten keywords,
linear systems, stochastic systems, L\'evy processes, system
identification, financial modelling              % chosen from the IFAC
\end{keyword}                             % keyword list or with the
                                          % help of the Automatica
                                          % keyword wizard
----
\end{comment}
\title{Identification of Finite Dimensional Linear Systems Driven by L\'evy processes}
\date{}
\author{L\'aszl\'o Gerencs\'er and M\'at\'e M\'anfay}

\maketitle
%----
%\begin{abstract}                          % Abstract of not more than 200 words.
L\'evy processes are widely used in financial mathematics,
telecommunication, economics, queueing theory and natural sciences
for modelling.
%Price processes are then defined as a corresponding geometric
%L\'evy process, implying the fact that returns are
%independent. In this paper we propose an alternative class of
%models allowing to describe dependence between return data.
A typical model is obtained by considering finite dimensional
linear stochastic SISO systems driven by a L\'evy process. In this
paper we consider a discrete-time version of this model driven by
the increments of a L\'evy process, such a system will be called
L\'evy system. We focus on the problem of identifying the dynamics
and the noise characteristics of such a L\'evy system. The special
feature of this problem is that the statistical description of the
noise is given by the characteristic function (c.f.) of the
driving noise not by its density function. As an alternative to
the maximum likelihood (ML) method we develop and analyze a novel
identification method by adapting the so-called empirical
characteristic function method (ECF) originally devised for
estimating parameters of c.f.-s from i.i.d. samples. Precise
characterization of the errors of these estimators will be given,
and their asymptotic covariance matrices will be obtained. We also
demonstrate that the arguments implying asymptotic efficiency for
the i.i.d. case can be adapted for the present case.
%----
%\end{abstract}

%----
%\end{frontmatter}

\section{Introduction}
\label{intro}

L\'evy processes are widely used to model phenomena arising in
natural sciences, economics, financial mathematics, queueing
theory and telecommunication
\cite{wind},\cite{novikov},\cite{CONT-TANKOV-FinModJunp}. The
classical model for modelling market dynamics, namely geometric
Brownian motion, was proposed by Louis  Bacehelier
\cite{Bacehelier}. This model is still the accepted core model
despite the fact that empirical studies revealed that its
assumptions are not realistic. For example, since price movements
are induced by transactions which can be unevenly distributed in
real time, it would be more natural to use a time changed Brownian
motion to model price dynamics. If the time change is defined by a
gamma process, we obtain the so-called VG (shorthand for Variance
Gamma) process. VG processes reproduce a number of stylized facts
of real price processes, such as fat tails and large kurtosis. It
can be shown that the above time changed Brownian process itself
is a L\'evy process. Extending the above construction novel price
dynamics have been proposed by a variety of authors, called the
geometric L\'evy processes obtained by exponentiating a L\'evy
process.
%%% ie. the Brownian motion is replaced by a L\'evy process.
%L\'evy processes have become a widely used tool in modeling price
%processes of financial instruments, such as stock prices or
%indices \cite{novikov}.

The objective of this paper is to present a combination of
advanced techniques in systems identification with a specific
statistical technique, widely used in the context in finance,
called the ECF (shorthand for empirical characteristic function)
method. The ECF method was originally designed for i.i.d. samples
and A.~Feuerverger and P.~McDunnogh \cite{FEUER+MC} showed that it
can be interpreted as the Fourier transform of an ML method.
Several papers study the problem of identifying the noise
characteristics of a linear system, but only a few pays attention
to the problem of identifying the system parameters as well.
Brockwell and Schlemm \cite{Brockwell-Schlemm} consider the
parametric estimation of the driving L\'evy process of a
multivariate continuous-time ARMA processes, but the
identification of system parameters is out of the scope of their
paper. Calder and Davis \cite{calder_davis} discuss the
M-estimators of ARMA processes with a given distribution on the
noise process. The quasi-maximum likelihood estimation of
multivariate L\'evy-driven continuous-time ARMA processes is
studied by Schlemm and Stelzer in \cite{quasi_max}, the method
presented there identify the system parameters and the covariance
structure of the noise process, but further characteristics of the
driving noise are not estimated.

In this paper we present a three-stage identification method for single-input-single-output (SISO)
that estimates both the system parameters and the characteristics of the noise process. We give the precise
characterization of the estimation error as well. We adapt the ECF method for linear systems and demonstrate
that our method can outperform standard system identification methods such as prediction error method or quasi
maximum likelihood estimation method and that it is asymptotically efficient.
%uj{\bf ADD REF BELOW !!!}
In \cite{ECC_own} the same problem is tackled. Two methods are proposed, the first
is a two-step method that combines the prediction error method and
the empirical characteristic function method for i.i.d. data, the
second estimates the the system parameters and the noise
parameters simultaneously. It is proved that the second method may
estimate the system parameters in a more efficient way than the
first one, still it does not give an efficient estimator.
Moreover, the method presented in
%uj{\bf ADD REF !!!}
\cite{ECC_own} is applicable only if the driving noise is a zero mean process.

%\medskip \noindent
%{\bf !!! BEGIN CORRECTIONS :}
%\medskip

A L\'evy process $(Z_t),~t \in \mathbb{R}$ is a continuous-time stochastic
process, which is much like a Wiener process: it is a stochastic
process with stationary an independent increments, but
discontinuities or
jumps are allowed. %A good survey paper on L\'evy processes used in
%financial modelling is the paper by Miyahara and Novikov,
%\cite{novikov}. \cite{raible} studies several problems arising in the field of exponential L\'evy processes.
For an excellent introduction to the theory of L\'evy processes
see \cite{Sato}.
%uj{\bf WHY NOT SATO ?}

A key building block in the theory of L\'evy processes is the
compound Poisson process, which is a Poisson process with random,
independent and identically distributed jumps. Extending the
notion of compound Poisson processes, a more general class L\'evy
processes is obtained via the formal expression
\begin{equation}
\label{eq:Levy_def}
Z_t = \int_0^t \int_{{\mathbf R}^1} x N(ds, dx),
\end{equation}
where $N(dt,dx)$ is a time-homogeneous, space-time Poisson point-
process, counting the number of jumps of size $x$ at time $t$. A
simple and elegant introduction to Poisson point-processes in a
general state-space is given in \cite{kingman}. The
interpretation of the above integral is fairly straightforward,
but attention should be paid to some technical conditions, see
below. The above process $(Z_t)$ is called a pure-jump L\'evy
process, indicating the fact that it is purely defined in terms of
its jumps.

A basic technical tool in the theory of Poisson point-processes is
the intensity of the process.
%
%
%%%In this case $Z_t$ is a pure jump process, which paradoxically
%%%means that the L\'evy-Ito decomposition of $Z_t$ does not have a
%%%Brownian motion component (but it may have a drift term).
%
%
%For an excellent survey on L\'evy process see {\bf cite{SATO}}, and
%{\bf cite{JACOD+SHIRYAEV}}.
%Then $Z_t$ can be written as
%$$
%Z_t = Z_t^0 + \gamma t,
%$$
%where $Z_t^0$ is a pure jump function:
%$$
%Z_t^0 = \sum_{0 \le s \le t} \Delta Z_s.
%$$
%
%
The intensity of a Poisson point-process $N(dt,dx)$ is formally
defined by $\E [N(dt,dx)],$ with $\E $ denoting expectation. Due
to time homogeneity, $\E [N(dt,dx)]$ can be written as
$$\E [N(dt,dx)] = dt \cdot \nu(dx),$$
where $\nu(dx)$ is the so-called L\'evy-measure.
%
%
%L\'evy processes are
%often best defined via their L\'evy-measures. In the case when
%$\int_{{\mathbf R}^1 \backslash 0} \nu(dx) < \infty$ we get a
%compound Poisson process with a finite number of jumps in each
%finite interval (or a finite activity process in the language of
%finance).
%
%
Now the above representation of a pure-jump L\'evy process given
in (\ref{eq:Levy_def}) is mathematically rigorous if
\begin{equation}
\label{eq:INTEGRABILIY}\int_{{\mathbf R}^1}  \min (|x|,1) \nu(dx)
< \infty. \end{equation} The intuition behind this condition that
small jumps have a finite contribution, and thus the fine
structure of the sample paths is relatively smooth. In fact, it
can be shown that the sample paths of $(Z_t)$ are of {\it finite
variation} with probability $1.$ In the area of financial time
series sample paths with finite variations are obtained for most
indices, as supported by empirical evidence, see \cite{CGMY}. This
phenomenon may be explained by the averaging effect when computing
an index, such as SP500. We note in passing that condition
(\ref{eq:INTEGRABILIY}) implies that for all finite $t$
\begin{equation}
\E \left[|Z_t|\right] < \infty.
\end{equation}
For an explanation see the comments following \ref{cond:MOMENTS_of_NU} below.

%uj\medskip \noindent {\bf !!! LABEL IT !  SAY,
%label{cond:MOMENTS-of-NU}}
%%%% or rather label{cond:$MOMENTS_of_NU$}}
%\medskip

It is easily seen that the characteristic function of $Z_t$ can be
written in the form
\[
%\label{eq:}
\E \left[e^{iuZ_t}\right] = e^{t \psi(u)}.
\]
Here  $\psi(u)$ is called the characteristic exponent. Note that
the logarithm of the characteristic function is linear in $t,$
which is implied by the fact that $(Z_t)$ has independent and
stationary increments.

The characteristic function plays a key role in the study of
L\'evy processes, because, unlike the density function of $Z_t,$
it typically has a closed form. The c.f. of a L\'evy process is
given by the celebrated L\'evy-Khintchine formula, which, in the
case of processes defined by (\ref{eq:Levy_def}), reduces to the
following:
\[
\E \left[ e^{iuZ_t} \right]=\exp \left[ t\left( ibu +
\int_{{\mathbf R}^1} \left( e^{iux}-1\right) \nu(dx) \right)
\right],
\]
where $b$ is known as the drift coefficient.
%uj{\bf ADD: EXPLAIN $b$ !!!}

%All other processes can be
%obtained by a limiting procedure in $L_2$, leading to the
%condition:
%$$ \int_{{\mathbf R}^1 \backslash
%0}  \min (|x|^2,1) \nu(dx) < \infty.
%$$

%\begin{lem} \cite{Sato}
%\label{prop:CHANGE_of MASURES_RN_PUREJUMP} Let $(X_t, P)$ and $(X_t, P')$
%be two L\'evy-processes on $\mathbf R$ with characteristic triplets
%$(A, \nu, \gamma_h)$ and with characteristic triplet $(A', \nu',
%\gamma_h')$, with respect to the same truncation function $h$. Assume that $P$ and
%$P'$ are equivalent on ${\mathcal F}_T$, for $T < \infty$, and the
%processes are diffusion-free. i.e. $A=A'=0$. Let
%$$ Y(x) =
%\zjel {\frac {d \nu'} {d \nu}}(x).
%$$
%Then the Radon-Nikodym derivative of $P'$ w.r.t. $P$ on ${\mathcal
%F_t}$ with $0 \le t < T$ is given by
%\begin{equation}
%\label{eq:}
%\left[ {\frac {dP'} {dP}}\right]_{\mathcal F_t}= e^{U_t}
%\end{equation}
%with
%\begin{equation}
%\label{eq:}
%U_t = \lim_{\epsilon \rightarrow 0} \left( {\sum_{s \le t, |\Delta X_s|
%> \epsilon} \log Y (\Delta X_s) - t \int_{|x|
%> \epsilon}
%(Y(x)-1) \nu(dx)} \right).
%\end{equation}
%\end{lem}

%!!! explain difficulty.

\subsection{Discrete-time L\'evy-systems}
\label{sec:disc_levy}
%The simplest example is a L\'evy AR(1) process,% defined by
%\begin{equation}
%dY_t=aY_{t-} dt+d Z_t
%\end{equation}
%with $a<0.$
% but it turns out that even in this simple case the
%identification problem is far from trivial. Therefore
%To avoid this problem we consider an alternative discrete-time
%model class.

A natural object for study is a linear stochastic system driven by
a L\'evy-process. Since the study of continuous-time systems would
lead to a number of technical difficulties, we restrict our
attention to discrete-time, finite-dimensional linear stochastic
systems driven by the increments of a L\'evy-process:
\begin{equation}
\label{eq:disc_levy2} \Delta y=A(\theta^*, q^{-1}) \Delta Z,
\end{equation}
%uj{\bf !!! A vs H ???}
%where $\Delta Z_n=(\Delta Z_n^{(1)},\cdot\cdot\cdot,\Delta Z_n^{(m)})$ is the increment of an $m$-dimensional L\'evy process $Z$ over
defined for the time range $- \infty < n < + \infty,$ where
$\Delta Z_n$ is the increment of a L\'evy process $(Z_t)$ with $
-\infty < t < + \infty$, and $Z_0=0$,
%
%\medskip \noindent
%uj{\bf !!! CHECK THE CONDITION FOR $Z$ ABOVE !
%Z_0=0 maradhat}
%\medskip
%
%
over an interval $[(n-1)h,nh),$ with $h>0$ being a fixed sampling
interval, and $ -\infty < n < + \infty$. The L\'evy-measure of $Z$
will be denoted by $\nu(dx)=\nu(dx,\eta^*),$ where $\eta^*$
denotes an unknown parameter-vector with a known range, say
$D_{\eta} \subset \mathbb{R}^r.$ The operator $A(\theta^*,
q^{-1})$ is a rational, stable (causal) and inverse stable
function of the backward shift operator $q^{-1},$ depending on
some unknown parameter-vector $\theta^*,$ taking its values from
some known set $G_{\theta} \subset \mathbb{R}^p.$
%\medskip \noindent
%uj{\bf THE NEXT CONDITION IS CERTAINLY INCORRECT, SINCE $A$ AND $H$
%ARE MIXED UP. WE MAY SIMPLY KEEP THE CONDITION AS PART OF THE
%INTRODUCTORY TEXT AND ADD DETAILS LATER FOR THE SS REP}
%\medskip
%\textbf{Condition 2}
% $A(\theta)$ is assumed to be exponentially stable and exponentially inverse stable
%for $\theta\in G_{\theta}\subset \mathbb{R}^p,$ where $G_{\theta}$
%is a known open set.
The observed output process is then $\Delta y.$ We call such
systems briefly L\'evy-systems.

The fundamental problem to be discussed in this paper is the
efficient identification of both the systems and the noise
parameters. The ML method would be appropriate in solving the full
identification problem (i.e. estimating both $\theta^*$ and
$\eta^*$) along standard lines, if we knew the density function of
$\Delta Z_n $ is known, see \cite{gerencs_ML} and
\cite{LSS_ML_GL+MGY+RZ}. Unfortunately, typically this is not the
case with the mostly used L\'evy processes.

Therefore we develop a new method, using a combination of the PE
(prediction error) and an adapted version of the so-called ECF
(empirical characteristic function) method, widely used in
finance, to get a competitive alternative to the ML (maximum
likelihood) method.

The ECF method was originally designed for i.i.d. samples. It has
the remarkable property that under certain idealistic assumptions
it is as efficient as the ML method. Certain extensions to
dependent data are available in the literature at the cost of
losing efficiency. Our main contribution is the development of a
method for system identification using a suitably adapted ECF
method, the efficiency of which is established solely relying on
efficiency results for i.i.d. data.

Let us now describe a few additional technical details of our
model. Let us assume that a state space representation in
innovation form equation for this model is given by
\begin{align}
\Delta X_{n+1}&=H(\theta^*) \Delta X_n + K(\theta^*) \Delta Z_n \\
\Delta Y_n&=L(\theta^*) \Delta X_n + \Delta Z_n.
\end{align}
Then stability and inverse stability of the system is then
described by the following condition:
\begin{cond} \label{cond:STAB} It is assumed that the system matrices $H(\theta)$
and $H(\theta) - K(\theta)L(\theta)$ are stable for $\theta \in
G_{\theta}.$
\end{cond}
To define the smooth dependence on $\theta$ suppose that $A(\theta, q^{-1})$ is
three-times continuously differentiable w.r.t. $\theta$ for $\theta
\in D_{\theta}$. Let $\mathscr{F}$ denote the natural filtration with $\mathscr{F}_{n-1}^{\Delta Z}=\sigma \left\{ \Delta Z_k: k \leq n-1 \right\}$.
The system (\ref{eq:disc_levy2}) is certainly well-defined if
$\Delta Z$ satisfies the integrability condition
(\ref{eq:INTEGRABILIY}). Namely, then $\Delta y_n$ can be written
as a weighted sum of past values of $\Delta Z,$ with exponentially
decaying weights, converging in $L_1.$
%\medskip \noindent
%uj{\bf !!! MOVE NEXT BIT DOWN TO ECF CHAPTER:}
%\medskip
%A nice heuristic justification for this has been given by
%A.~Feuerverger and P.~McDunnogh in \cite{FEUER+MC}, showing that
%the equations defining the ECF method for i.i.d samples can be
%obtained as the Fourier transform of the likelihood equations.
%
%
%%%The reason behind that we work with discrete time systems is that
%%%in the statistical analysis of continuous time systems the
%%%difficulty of applying a maximum-likelihood (ML) method lies in
%%%the fact that there is no natural reference measure in the space
%%%of sample paths. In addition, the computation of the Radon-Nikodym
%%%derivative is practically not feasible since $\int_{-\infty}^{t}
%%%e^{H(t-s)} dZ_s $ is not even a L\'evy process.
%
%
%With the same conditions that will be discussed in Section VI. the following result holds.
%\medskip \noindent
%uj{\bf USE LABELS !!!}
%\medskip
We will need the following technical condition:
%\medskip \noindent
%uj{\bf HOW DO WE DEAL WITH THE CASE $Q=\infty$ ???. Shall we write
%$q < Q$ ALL THE WAY ??? JUSTIFICATION FOR $Q$ ???}
%\medskip
\begin{cond}\label{cond:MOMENTS_of_NU}  We assume that for  all $q \geq 1$
  \begin{equation}
  \label{eq:moments}
  \int_{|x| \ge 1}|x|^q\nu(dx)<+\infty,
  \end{equation}
  and that $\E \left[ \Delta Z_n\right]=0.$
\end{cond}
It follows, see \cite{Sato}, that for $q \geq 1$ and for all $h \ge 0,$ the $q$-th
moments of the increments of $Z$ are finite:
\begin{equation}
\E \left[|\Delta Z_n|^q \right]< \infty.
\end{equation}
%
%uj{\bf IF WE HAVE $t \ge 0$ IN DEFINING A LEVY PROCESS, THEN THE
%LEVY-SYSTEMS SHOULD ALSO BE DEFINED ONLY FOR $t \ge 0$ !!!}
%
%
We note here that Condition \ref{cond:MOMENTS_of_NU} holds in
our benchmark examples to be presented in the next Section.
%
% The standard examples we can think of is that $Z$ is
%  an $\alpha$-stable process with
%
%$$\nu(dx)=\frac A {x^\alpha}\mathbf{1}_{x>0}dx+\frac B
%  {|x|^\alpha}\mathbf{1}_{x<0}dx,$$
%
%or $Z$ is a tempered stable process, such as CGMY with
%
%$$\nu(dx)=\frac{Ce^{-G|x|}}{|x|^{1+Y}}\mathbf{1}_{x<0} +
%\frac{Ce^{-M|x|}}{x^{1+Y}}\mathbf{1}_{x>0}.$$
%
%{\bf SOPPED HERE. SAT. 20:20/. ssssssssssssssss}
%%%The transfer function $A(\theta)$ is the rational function of the backward shift operator:
%%%$$
%%%A(\theta,z^{-1})=\sum_{i=0}^{\infty} a_{i}(\theta) z^{-i},
%%%$$
%%%where $z^{-1}\Delta Z_n=\Delta Z_{n-1}.$

Let $D_{\theta}$ and $D_{\theta}^*$ be compact domains such that
$$
\theta^* \in D_{\theta}^* \subset \text{int } D_{\theta} \quad
{\rm and} \quad D_{\theta} \subset G_{\theta}.$$ The domains
$D_{\eta},D_{\eta}^*$ are defined analogously. Finally, let
$\rho^*=(\theta^*,\eta^*)$ denote the joint parameter vector, and
set
$$
D_{\rho} = D_{\theta} \times D_{\eta}, \qquad D_{\rho}^* =
D_{\theta}^* \times D_{\eta}^*,  \qquad G_{\rho} = G_{\theta}
\times G_{\eta}.
$$
%
%
%%%\textbf{Condition 2}
%%% $A(\theta)$ is assumed to be exponentially stable and exponentially inverse stable
%%%for $\theta\in G_{\theta}\subset \mathbb{R}^p,$ where $G_{\theta}$ is a known open
%%%set.
%
%
%%A system is exponentially stable if all the eigenvalues of $A$
%%have strictly negative real parts.
%%%
%
%
%
Before going into further details we present a few examples of
L\'evy processes used for modeling purposes.
%
%
\begin{comment}
If we knew the probability density function of the noise $\Delta
Z_n$ then we could apply an ML (Maximum Likelihood) estimation
method, and establish sharp results for the estimation error, see
\cite{LSS_ML_GL+MGY+RZ}. The challenge of the present problem is
that it is the characteristic function of the noise that is
explicitly given. A natural approach to solve this problem is to
combine techniques of system identification with the empirical
characteristic function (ECF) method widely used in finance to
analyze i.i.d. data. Before going into further details we present
a few examples of L\'evy processes used for modeling purposes.
\end{comment}
%In the statistical analysis of such systems, both the system
%dynamics and the fine characteristics of $(Z_t)$ are to
%identified. The difficulty of applying a maximum-likelihood (ML)
%method lies in the fact that there is no natural reference measure
%in the space of sample paths. To avoid this problem we focus on a
%simpler class of discrete time finite dimensional L\'evy systems.
%as opposed to the mainstream of the system identification literature.
%The objective of this paper is to estimate the dynamics and noise
%characteristic of $Z_t$ using historical data. Obviously,
%observing the price $S_t$, we also have observed $Y_t.$ Thus our
%problem boils down to estimating the dynamics and noise
%characteristic of a L\'evy system.
\section{Examples for widely used L\'evy processes}
\label{sec:levy_examples}
%L\'evy-measures nicely reflect a number of important properties of
%properties of $Z_t= Z_t^0$ itself. Thus e.g. $$\E e^{\lambda X_t}
%< \infty$$ for all $t >0 $ or equivalently for some $t >0$ if and
%only if
%$$
%\int_{|x|\ge 1} e^{\lambda x} \nu (dx) < \infty.
%$$
%In this case $\E e^{\lambda X_t}  = e^{t \psi(\lambda/i)} = e^{t
%\psi(-i \lambda)}$, see Theorem 25.17 in Sato's book: cite{SATO}
%%% Get Sato!

%uj{\bf ADD EXAMPLES FROM TELECOMMUNICATION (TIEN) OR ROAD SURFACE MODELLING ???}
Compound Poisson process is defined by a rate $\lambda$ and a jump size distribution $F$ via
$$
Z_t=\sum_{i=1}^{N_t}X_i,
$$
where $N_t$ is a Poisson process with rate $\lambda,$ and $X_i$-s
are i.i.d. random variables with distribution $F.$ Compound
Poisson processes are widely used for modeling in queueing theory.
For example in \cite{tien} a generalized multi-server queue is
used to model telecommunication networks. Among several properties
of the model, customer arrivals, server failures and packet losses
are modeled with compound Poisson processes.

To model the increments of the logarithm of a price process a wide
range of geometric L\'evy processes has been proposed by a variety
of authors.
%Particular attention has been paid to price
%dynamics driven by so-called VG (shorthand for Variance Gamma)
%processes, proposed by Madan, Carr and Chang \cite{CARR-VG}.  A VG-process
%is a time changed Brownian motion with drift, where the time
%change is a so-called gamma process, which is essentially the
%continuous time extension of the inverse of a Poisson process.% It
%is a three-parameter class of processes, with explicit
%characteristic function and L\'evy measure.
%The logic of using a
%time changed Brownian motion, as an alternative, is that price
%movements are induced by actual transactions which can be very
%uneven in real time. Geometric VG processes reproduce a number of
%so-called stylized facts of real price processes, such as fat
%tails or high kurtosis.
%
Mandelbrot suggested to use $\alpha$-stable process to model the
price dynamics of wool, see \cite{MANDELBROT}. An $\alpha$-stable
with $0 < \alpha < 2$ is defined via the L\'evy measure
\begin{equation}
\nonumber
% \label{eq:}
\nu(dx) = C^- |x|^{-1-\alpha} {\mathbf 1}_{x < 0 } dx+ C^+ |x|^{-1-\alpha} {\mathbf 1}_{x
> 0 } dx.
\end{equation}
%Note that to ensure the integrability of $x^2 \nu(dx)$ around zero
%we must have $\alpha <2.$
%
A recently widely studied class of L\'evy processes is the CGMY
process due to Carr, Geman, Madan and Yor \cite{CGMY}. It
is obtained by setting $C^- = C^+$, and then, separately for $x
>0$ and $x<0$, multiplying the
L\'evy-density of the original symmetric stable process with a
decreasing exponential. The corresponding L\'evy-measure, using
standard parametrization, is of the form:
\begin{equation*}
\nonumber
%\label{eq:}
\nu(dx) = {\frac {C e^{-G |x|}} {|x|^{1+Y}}}  {\mathbf 1}_{x < 0 } dx +
{\frac {C e^{-M x}} {|x|^{1+Y}}} {\mathbf 1}_{x
> 0 } dx,
\end{equation*}
where $C,G,M >0$, and $0 < Y <2$. Intuitively, $C$ controls the
level of activity, $G$ and $M$ together control skewness.
Typically $G > M$ reflecting the fact that prices tend to increase
rather than decrease. $Y$ controls the density of small jumps,
i.e. the fine structure. For $Y<1$ the integrability condition (\ref{eq:INTEGRABILIY}) is satisfied, thus corresponding L\'evy
process is of finite variation. The characteristic exponent $\psi(u)$ of the
CGMY process is given by
\begin{equation*}
%\label{eq:}
C \Gamma (-Y) \left((M-iu)^Y- M^Y + (G+iu)^Y- G^Y\right),
\end{equation*}
where $\Gamma$ denotes the gamma-function. A more general class of tempered stable distributions is studied in $\cite{terdik},$ Terdik and Woycz\'nski obtains analytic formulas for the Rosi\'nski measure of tempered processes.

%Allowing $C$ and $Y$ to take on different values for $x
%>0$ and $x<0$ we get a more general class of processes called tempered stable
%process. see cite{}.

%uj VG added
Formally setting $Y=0$ we get the L\'evy density of the so-called
Variance Gamma process (VG for short) that has been proposed by
Madan, Carr and Chang \cite{CARR-VG}. The Variance Gamma (VG)-process
is a time changed Brownian motion with drift, where the time
change is a so-called gamma process, which is essentially the
continuous time extension of the inverse of a Poisson process. It
is a three-parameter class of processes, with explicit
characteristic function and L\'evy measure.
Let
$B_t(\theta,\sigma)$ be a Brownian motion with drift $\theta$ and volatility $\sigma$, i.e.:
$$B_t(\theta,\sigma)=\theta t+\sigma B_t,$$ where $(B_t)$ is a standard Brownian motion and let
$\gamma_t(\mu,\nu)$ be a gamma process with mean $\mu$ and
variance $\nu,$ i.e. $\gamma_t(\mu,\nu)$ is a stochastic process
with independent gamma distributed increments.
%uj2 explained via skewnes, kurtosis
The VG process $\left(X_t(\sigma,\nu,\theta)\right)$ is defined as
$$X_t(\sigma,\nu,\theta)=B_{\gamma_t(1,\nu)}(\theta,\sigma).$$
Hence, the VG process is a time-changed Brownian motion. According to \cite{CARR-VG} $\theta$ controls the skewness and $\nu$ controls the kurtosis of the process. A possible other definition of the VG process is that $X_t(\sigma,\nu,\theta)$ can be written as the difference of two gamma processes.

The L\'evy measure of a VG process can be obtained by first
computing its characteristic function and then applying
L\'evy-Khintchine's formula in the inverse direction. Thus we get:
\[
%\nonumber
%\label{eq:}
\nu(dx) = \left\{
        \begin{array}{ll}
                \frac{\mu^2_{n}}{\nu_n} \frac{\exp ( -\frac{\mu_n}{\nu_n}|x|)}{|x|}dx  & \mbox{if } x < 0 \\
                \frac{\mu^2_{p}}{\nu_p}\frac{\exp (-\frac{\mu_p}{\nu_p}x)}{x}dx  & \mbox{if } x > 0,
        \end{array}
\right.
\]
where the parameters $\mu_p, \nu_p, \mu_n,\nu_n$ are obtained in
terms of the original parameters as follows:
\begin{eqnarray*}
\mu_p=\frac12 \sqrt{\theta^2+\frac{2\sigma^2}{\nu}}+\frac{\theta}{2} \quad \nu_p=\mu^2_p \nu \\
\mu_n=\frac12 \sqrt{\theta^2+\frac{2\sigma^2}{\nu}}-\frac{\theta}{2} \quad \nu_n=\mu^2_n \nu
\end{eqnarray*}
From here we get the following remarkable property of VG
processes: a VG process $X_t(\sigma,\nu,\theta)$ can be written as
the difference of two gamma processes $(\gamma_{p,t})$ and $(\gamma_{n,t})$ :
\[
X_t(\sigma,\nu,\theta) = \gamma_{p,t}(\mu_p,\nu_p)  -
\gamma_{n,t}(\mu_n,\nu_n).
\]
In particular, it follows that a VG process is of finite
variation.

%uj{\bf I WANT MORE FLESH !}
%TIEN added, more VG stuff
\section{The empirical characteristic function method for i.i.d. data}

We briefly describe the ECF method for i.i.d.
samples, see \cite{CARRASCO-EFFEMPIRCHAR}. % It is worth mentioning that the method can handle
%Markov-dependent processes as well, see \cite{CARRASCO-EFFCGMM}.
The ECF method gives an efficient estimate of the unknown parameters of a given family of distributions \cite{CARRASCO-EFFEMPIRCHAR}.
%uj moved here!
A nice heuristic justification for this has been given by
A.~Feuerverger and P.~McDunnogh in \cite{FEUER+MC}, showing that
the equations defining the ECF method for i.i.d samples can be
obtained as the Fourier transform of the likelihood equations.

Let $(r_1, r_2, \ldots r_N)$ be i.i.d. observations, and let a
closed form of the characteristic function $\varphi(u, \eta)$ be
known, with $\eta$ being a $p$-dimensional parameter vector, and $u \in
\mathbb{R}$. The true value of the parameter will be denoted by
$\eta^*$.

The idea is to estimate $\eta^*$ by a value of $\eta$ for
which the characteristic function (cf) best matches the empirical
characteristic function (ecf). The error for any fixed $u$ is
defined as
\[
\overline{h}_N(u, \eta)=\frac 1N \sum_{k=1}^N h_k(u, \eta),
\]
where $h_k(u, \theta)$ is the generalized moment function:
\[
h_k(u, \eta) = e^{iu r_k}-\varphi(u, \eta).
\]
An important property of the moment function is that $$\E
\left[h_k(u, \eta^*)\right]=0, \quad {\rm for~all} \quad u,$$ where
$\eta^*$ denotes the true parameter.
In the case of a finite set of moments evaluated, say at $u_1,
..., u_M$, with $M > p$, define
$$\overline{h}_N(\eta)= (\overline{h}_N(u_1, \eta),...,\overline{h}_N(u_M, \eta))^T,$$ and its expectation
$$
g(\eta) = \E \left[\overline{h}_{N}(\eta)\right].
%%% = \varphi(u_k,\eta^*)-\varphi(u_k,\eta).}
$$
Since $g(\eta^*)=0,$
we would obtain $\eta^*$ by
solving the over-determined system of equations
\begin{equation}
g(\eta)=0,
\end{equation}
where ${\rm dim~ } g(\eta)=M>p.$
Since $g$ is not computable in practice we seek a least-square solution by minimizing the weighted quadratic cumulative error
\begin{equation}
\label{eqn:Q_FINITE} V_N(\eta)=|K^{-\frac12} \overline{h}_N(\eta)|^2
\end{equation}
where $K$ is an appropriate, $m \times m$ weighing matrix.

Now we compute the asymptotic covariance matrix of the estimated parameter $\hat{\eta}_N.$
We denote the complex conjugate of a vector or matrix $M$ with $M^*.$
The gradient equation ($p$ equations):
$$
{\overline h}^*_{\eta,N} (\eta) K^{-1} {\overline h}_N (\eta) = 0 .
$$
The corresponding approximating problem can be defined via
$$
G K^{-1} {\overline h}_N (\eta) = 0,
$$
where $G$ is the $M \times p$ matrix
$
G = g_{\eta} (\eta^*).
$
The left hand side is considered as a new set of exactly $p$ scores.
The asymptotic gradient is
$$
g^*_{\eta} (\eta) K^{-1} g (\eta),
$$
while its derivative at $\eta^*$ (the Hessian of the asymp. cost) is
$$
R= g^*_{\eta} (\eta^*) K^{-1} g_{\eta} (\eta^*).
$$
Then the Hessian of the asymp. cost is
$$
T= G^* K^{-1} G.
$$
To get the asymptotic covariance of the new set of scores define the $M \times M$ covariance matrix by
$$
C_{k,l} = \E \left[ ~h_{n}(u_k,\eta^*) h^*_{n}(u_l,\eta^*) \right].
$$
Note that we have
$$
C_{k,l} = \varphi(u_k - u_l,\eta^*) - \varphi(u_k ,\eta^*) \varphi(- u_l,\eta^*).
$$
Thus the asymptotic covariance of the new set of scores is
$$
S = G^* K^{-1} C K^{-1} G .
$$
The asymptotic covariance of the estimator $\hat \eta_N$ is then
$$
T^{-1} S T^{-1},
$$
or equivalently,
$$
(G^* K^{-1} G)^{-1} ~G^* K^{-1} C K^{-1} G ~(G^* K^{-1} G)^{-1}.
$$
It is easy to see that the optimal value of $K$ is
$$
K= C
$$
yielding the asymptotic covariance for $\hat \eta_N$
$$
(G^* C^{-1} G)^{-1}.
$$

%\section{The identification problem}
%\label{sec:problems} In this section we re-formulate the
%identification problem for discrete-time, finite dimensional
%L\'evy systems, and formulate two related, simpler identification
%problems which are of technical interest. We also sketch a the
%proposed path to their solutions.
%
%%\textit{Both the system parameters and the noise parameters are
%%unknown.}
%
%The main problem that we consider is the identification of L\'evy
%systems, when both the system parameters and the noise parameters
%are unknown. In Section \ref{sec:discussion} we briefly sketch a method proposed in the
%literature, see \cite{knight}, \cite{YU-EMPIRCHARFUN},
%%
%%\medskip
%%\noindent {\bf !!! DO PROVIDE SOME REFERENCES !}
%%\medskip
%%
%for the statistical analysis of dependent data via the ECF method,
%and also point out the shortcoming of this approach.
%

%
%
%\medskip \noindent
%{\bf  !!! PROVIDE EXACT REF !}
%\medskip

\section{Three-stage method}\label{stage1-3}

The main problem that we consider is the identification of L\'evy
systems, when both the system parameters and the noise parameters
are unknown. In Section \ref{sec:discussion} we briefly sketch a method proposed in the
literature, see \cite{knight}, \cite{YU-EMPIRCHARFUN},
%
%\medskip
%uj\noindent {\bf !!! DO PROVIDE SOME REFERENCES !}
%\medskip
%
for the statistical analysis of dependent data via the ECF method,
and also point out the shortcoming of this approach.

Therefore we propose a completely different approach, combining
the PE method and an adapted version of the ECF method. Our novel
method a three-stage method, the first stage being a standard PE
method for estimating the system parameters, taking into account
only that the innovation process is i.i.d. with having moments of
appropriate orders. Thus we get an estimate of the system
parameters, say, $\hat{\theta}_N.$ %(together with an estimate of
%the mean of the innovation).
%
%\medskip \noindent
%uj{\bf  !!! PROVIDE A NOTE ON NONZERO EXPECTATION !}
%\medskip
%we have zero expectation now

In the second stage, using a certainty equivalence argument,
pretending that $\hat{\theta}_N = \theta^*,$ we estimate the
innovation process by inverting the system using the estimated
system parameters $\hat{\theta}_N$. Then, the noise parameters
$\eta^*$ are estimated using the ECF method for i.i.d. sequences,
resulting in an estimate $\hat{\eta}_N.$ These procedures will be
briefly described in Section \ref{sec:STAGE-1-2}.
%
%\medskip \noindent
%{\bf  !!! CREATE SECTION STAGE-1-2 !}
%\medskip
%
Finally, in the third stage, once again using a certainty
equivalence argument, pretending that $\hat{\eta}_N = \eta^*,$ we
re-estimate the system parameters using a specific adaptation of
the ECF method for systems-identification with i.i.d. innovation
process, having a {\it known} characteristic function. This is the
most original step in our procedure.

The analysis of the effects of the estimation errors of
$\hat{\theta}_N$ and $\hat{\eta}_N$ on subsequent steps are based
on moment estimates of the estimation errors. The latter can be
obtained by extending the techniques of \cite{gerencs_arma}, and
exploiting the fact that all finite moments of the innovation
process are finite.
%This method will be studied in
%Section \ref{sec:mixed}.
%%%This procedure will be analyzed in Section \ref{sec:mixed}.
%%%The second, simplest problem is
%%%seemingly of
%%%mere technical interest:
%%%
%%%\textit{Known system parameters, unknown noise parameters.} In
%%%this case define and compute
%%%$$\varepsilon_n (\theta^*)=A^{-1}(\theta^*)\Delta y_n=A^{-1}(\theta^*)A(\theta^*)\Delta Z_n=\Delta
%%%Z_n,$$ assuming, for the sake of simplicity, that $\Delta Z_n =
%%%\varepsilon_n (\theta^*) = 0$ for $n \le 0.$ After that we can
%%%apply the ECF method for i.i.d. samples to obtain the estimation
%%%of $\eta^*.$ This simple solution will be the base of the second step of the identification method presented %in Section \ref{sec:mixed}.
%%%%%%see for example \cite{CARRASCO-EFFEMPIRCHAR}.

%uj{\bf !!! MATE: THIS REFERENCE SHOULD APPEAR ALSO IN THE INTRO !}

To set the stage for the final step of our procedure we briefly
summarize a simple known result on the ML estimate for the
identification of a linear stochastic system with i.i.d.
innovation of known characteristics, more accurately of known
probability density function, say $f(., \eta^*)$, following
\cite{gerencs_ML}. In this case we can obtain the maximum
likelihood estimate of the unknown system parameters $\theta^*$
via solving
\begin{equation}
\sum_{n=1}^{N} {\frac {\partial} {\partial \theta}} \log f
\left(\varepsilon_n(\theta),\eta^*\right)=0,
\end{equation}
where
\begin{equation}\label{eq:inverse}
\varepsilon_n(\theta)=A^{-1}(\theta)\Delta y_n
\end{equation}
is the estimated innovation process of the SISO system given under (\ref{eq:disc_levy2}).
%\medskip \noindent
%uj{\bf  !!! OF THE THE THE SISO SYSTEM GIVEN UNDER ... }
%\medskip
%
%
%\medskip \noindent
%{\bf  !!! COPY UP "SISO" TO INTRO ! }
%\medskip
%
%{\bf MATE: $\varepsilon_n(\theta)$ YET UNDEFINED !}

Then under certain technical conditions, such as the condition
that $\E[\Delta Z_n]= 0,$ the asymptotic covariance matrix of the
ML estimate
%%% $\hat{\theta}_N$
is given by
\begin{equation}
\label{eq:mu} \Sigma_{\rm ML}=\mu^{-1} \left(R_P^*\right)^{-1},
\end{equation}
where
$$
\mu=%\lim_{n \rightarrow \infty}
\E \left[\left(\frac{f'(\Delta Z_n,\eta^*)}{f(\Delta
Z_n,\eta^*)}\right)^2\right],$$
with $f'$ being the derivative of
$f(.,\eta^*)$ w.r.t its first variable, and $R_P^*$ is defined in connection with the PE method. Since this method is efficient we have that $\mu^{-1} \le \sigma^2,$ we note that the accuracy of the ML method can significantly
surpass that of the PE method, i.e. we can have $\mu^{-1}<<
\sigma^2.$ Large difference between $\mu^{-1}$ and $\sigma^2$ can be achieved by taking the mixture of a mass like continuous pdf with and another continuous pdf.
This $\mu$ can be interpreted as a Fisher of location parameter estimate.
This property of $\mu$ will be discussed in more details in Section \ref{sec:efficiency} and will be used in
analyzing the efficiency of our proposed method in Section \ref{sec:ECF_system}.
%\medskip \noindent
%{\bf  !!! $R^*$ AND $\varepsilon$ MAY BE DEFINED ABOVE IN
%CONNECTION WITH PE. "THIS IS WHAT PE WE CAN ACHIEVE." MAY REMARK
%THAT $\mu^{-1} \le \sigma^2$ ! }
%\medskip
%$$
%R^*=\lim_{n \rightarrow \infty} %\E\left[\varepsilon_{n\theta}(\theta^*)\varepsilon^T_{n\theta}(\theta^*)\right].
%$$
%
%\medskip \noindent
%{\bf  !!! ADD INTERPRETATION OF $\mu$ AS FISHER OF LOCATION. }
%\medskip
The challenge we address in this paper if we can achieve the same
accuracy in estimating $\theta^*$ when we know the characteristic
function of the innovation only rather than its p.d.f. The
surprising answer is a yes, or rather an almost yes. Before going
into details  we briefly summarize the first two stages of our
algorithm.

Hence, our proposed three-stage method can be summarized as follows
\begin{enumerate}
  \item Estimate $\theta^*$ by applying PE method, obtain $\hat{\theta}_N$
  \item Invert the system with $\theta=\hat{\theta}_N,$ then estimate $\eta^*$ by using the idea of ECF method for i.i.d. data and obtain $\hat{\eta}_N$
  \item Re-estimate $\theta^*$ by applying the ECF method for system identification to obtain an efficient estimate $\hat{\hat{\theta}}_N$ for the dynamics.
\end{enumerate}
%%%In our case, the p.d.f. of the noise distribution is not known.
%%%One might apply the prediction error method to estimate the system
%%%dynamics, i.e. $\theta^*.$ However, we will show that we may
%%%estimate $\theta^*$ in a more efficient
%%%way using an appropriate adaptation of the ECF method. This method will be presented and analyzed in Section %\ref{sec:ECF}. %In fact,
%this result is a special case of a more general result obtained
%for the general problem to be described in the next subsection.
%\section{Three-stage method: combining PE and ECF estimators}
%\label{sec:mixed}
%\label{STAGE-1-2}
\section{A summary of results on the PE method and the ECF method}
\label{sec:STAGE-1-2} %\label{sec:STAGE-1}

%uj{\bf  !!! PROVIDE A NOTE ON NONZERO EXPECTATION ! SEE MY SLIDES
%FOR THE ZABCZYK CONF.}
%we have zero exp. now, see note in 2. sentence.

%We identify $\theta^*$
%using only the orthogonality of $\Delta Z$ by applying a
%prediction error method. This way we get an estimation
%$\hat{\theta}_{N}$  of $\theta^*$, without using the characteristic function of $\Delta Z.$ Then we apply an ECF method with the score
%function
%$$h_n(u,\eta)=e^{iu\varepsilon_n(\hat{\theta}_{N})}-\varphi(u,\eta)$$
%to estimate $\eta^*.$

In this section we briefly summarize the PE method for the case
when the input noise has zero expectation, i.e. $\E[\Delta Z_n]=
0,$ and present a result that will be relevant later. Although
general L\'evy processes presented in Section
\ref{sec:levy_examples} are not zero mean processes, by
preprocessing our data, as is customary in classic time series
analysis, we may achieve that $\E[\Delta Z_n]= 0.$ First, we
define the estimated innovation process $\varepsilon(\theta)$ as
in the previous sections. The prediction error estimator of
parameters vector $\theta^*$ is then obtained by minimizing the
cost function
\begin{equation*}
V_{P,N}(\theta)=\frac{1}{2}\sum_{n=1}^N
\varepsilon^2_n(\theta),
\end{equation*}
over $G_{\theta},$
%uj2 changed
%{\bf ??? $G_{\theta} R,$ ???}
%
see \cite{ljung}. An alternative, more convenient definition of
the PE estimator $\hat{\theta}_N$ is obtained by setting the
gradient of the cost functions equal to zero, and considering the
equations:
\begin{equation*}
{\frac {\partial}  {\partial \theta}} V_{P,N}(\theta,
m)=\sum_{n=1}^N \varepsilon_{n\theta}(\theta)
~\varepsilon_n(\theta) =0
\end{equation*}
The asymptotic cost function associated with the PE method is defined as
$$
W_{P}(\theta)
%%% =\frac{1}{2}\lim_{n \rightarrow \infty} \E \varepsilon_n^2(\theta)
=\frac{1}{2} \E \left[ \left(\varepsilon^{(s)}_n(\theta)\right)^2 \right],
$$
where $\varepsilon^{(s)}_n(\theta)$ denotes the stationary
solution of (\ref{eq:inverse}) when $-\infty<n<\infty.$ (In
general, the superscript $^{(s)}$ will be used throughout this
paper if the marked stochastic process is obtained by passing
through a stationary process through an exponentially stable
linear filter starting at $- \infty$, as opposed to initializing
the filter at time $0$ with some arbitrary initial condition,
which is typically zero).
%%In general, the
%notation $(.)^*$ will be used throughout this paper if the
%corresponding stochastic process is obtained by passing through a
%stationary process through an exponentially stable linear filter
%starting at $- \infty$, as opposed to initializing the filter at
%time $0$ with some arbitrary initial condition, which is typically
%zero.
We have
\begin{equation*}
\begin{split}
&{\frac {\partial}  {\partial \theta}} W_{P}(\theta^*)=0 \text{ and } \\
&R_P^*:=W_{P,\theta\theta}(\theta^*)=\E \left[\varepsilon^{(s)}_{n
\theta}(\theta^*)\varepsilon^{(s)T}_{n \theta}(\theta^*)\right].
\end{split}
\end{equation*}
%uj {\bf !!! FIX $R^*$ and $\Sigma_{P}$!}
%{\bf !!! MATE: DID WE EXPLAIN WHAT IS $^*$ FOR?}
Furthermore, $\theta=\theta^*$ is the unique solution of $W_{P,\theta}(\theta)=0$ in $D_{\theta},$ see \cite{soderstrom}. The asymptotic covariance matrix of the PE estimate of $\theta^*$
is given by
%uj{\bf !!! STATE A THEOREM BELOW ! BEST OPTION: STATE THE MARTINGALE
%REP THM !}
\begin{equation}\label{eq:sigmap1}
\Sigma_{P}=\sigma^2\left(\E\left[\varepsilon^{(s)}_{n\theta}(\theta^*)\varepsilon^{{(s)}T}_{n\theta}(\theta^*)\right]\right)^{-1},
\end{equation}
where $\sigma^2$ is the variance of $\Delta Z_1.$
We will use this notation in a more general way:
\begin{defn}
For a stochastic process $X_n,$ and a function
$f:\mathbb{Z}\rightarrow \mathbb{R}^+$ we say that
$$ X_n=O_{M}(f(n)) $$ if for all $q \geq 1$
$$ \sup_{n} \frac{\E^{1/q}\left|X_n\right|^q}{f(n)} < \infty $$
holds.
\end{defn}
%uj{\bf !!! THIS WAS A STATEMENT ABOVE, NOW IT IS A CONDITION ??? }
%no condition
For the definition of $L$-mixing processes and for other corresponding definitions and theorems see the Appendix.
Theorem \ref{thm:theta_1_diference}, with minor variation, can be found in
\cite{gerencs_arma}.
%{\bf !!! MATE: EXPLAIN NOTATIONS BELOW!}
\begin{thm}\label{thm:theta_1_diference}
Under Conditions \ref{cond:STAB}, \ref{cond:MOMENTS_of_NU} we have
$$
\hat{\theta}_N-\theta^*=-(R_P^*)^{-1}V_{P,N,\theta}(\theta^*,m^*)+r_N,
$$
with $r_N=O_M(N^{-1}).$
\end{thm}
Direct consequence of this theorem is the following lemma, which can be proved with different methods as well.
\begin{lem}\label{lemmma:PE_N^-1/2}
Under Conditions \ref{cond:STAB}, \ref{cond:MOMENTS_of_NU} we have
$$
\hat{\theta}_N-\theta^*=O_M(N^{-1/2}).
$$
\end{lem}
%For the definition of $O^{Q/2}_{M},$ and for some details on $L$-mixing processes see the Appendix.

%uj{\bf ECF FOR I.I.D. DATA: I THINK THIS SHOULD BE DISCUSSED
%INDEPENDENTLY FOR DIDACTIC REASONS. RIGHT AFTER THE LEVY SECTION !
%THE TEXT BELOW THEN CAN BE ABRIDGED. }

\section{The ECF method for i.i.d. data. Application for estimating the noise parameters.}
\label{sec:ECF_iid_application}

The second, simplest problem is seemingly of mere technical
interest, when we know the system parameters, but the noise
parameters are unknown. In this case define and compute
$$\varepsilon_n (\theta^*)=A^{-1}(\theta^*)\Delta y_n=A^{-1}(\theta^*)A(\theta^*)\Delta Z_n=\Delta
Z_n,$$ assuming, for the sake of simplicity, that $\Delta Z_n =
\varepsilon_n (\theta^*) = 0$ for $n \le 0.$ After that we can
apply the ECF method for i.i.d. samples to obtain the estimation
of $\eta^*.$ An ideal score function for the ECF method to
estimate $\eta^*$ would be defined by
%uj2 i have changed again
%{\bf I have changed notations below: }
%
\begin{align}
h_{k,n}^{opt} (\theta^*,\eta) = e^{iu_k
\varepsilon^{(s)}_n(\theta^*)}-\varphi(u_k,\eta).
%h^*_n(u,\theta,\eta)&=e^{iu\varepsilon^*_n(\theta)}-\varphi(u,\eta)
\end{align}
Since we are not given $\theta^*$ we define an alternative,
$\theta$-dependent score function via
\begin{equation*}
h_{k,n}(\theta,\eta)=e^{iu_k\varepsilon_n(\theta)}-\varphi(u_k,\eta),
%h^*_n(u,\theta,\eta)&=e^{iu\varepsilon^*_n(\hat{\theta}_N)}-\varphi(u,\eta).
\end{equation*}
with a fix set of real numbers $u_1,\cdots,u_M$, with $NM \ge {\rm
dim}~\eta$. These are appropriate score functions since $$\E
\left[h^{(s)}_{k, n}(\theta^*,\eta^*)\right] =0.$$ Define
$$h_n(\theta,\eta)=\left(h_{1,n}(\theta,\eta),\cdots,h_{M,n}(\theta,\eta)\right)^T$$ and
$$\overline{h}_N(\theta,\eta)=\sum_{n=1}^N h_n(\theta,\eta).$$
The expectation of the score vector is denoted by
$$g(\theta,\eta)=\E\left[ \overline{h}_N(\theta,\eta)\right].$$
For a fixed $\theta$ we proceed like in the i.i.d. case and obtain
the $\theta$-dependent estimate $\hat{\eta}_N(\theta)$ of $\eta^*$
by finding a least squares solution to the over-determined system
of equations
$$
g(\theta,\eta)=\E\left[\overline{h}_N(\theta,\eta)\right]=0.
$$
More precisely, define the $\theta$-dependent cost function
$$
V_{E,N}(\theta,\eta)=\left|K^{-1}\overline{h}_N(\theta,\eta)\right|^2,
$$
%uj{\bf THE CLARIFICATION AND DETAILS I GAVE IN THE ZABCZYK TALK ARE
%MISSING, SUCH AS $G$ AND $C$ AND OPTIMAL $K$ ! }
where $K$ is a symmetric, positive definite weighting matrix and obtain $\hat{\eta}_N(\theta)$ as the solution of
$$
V_{E,N,\eta}(\theta,\eta)=0.
$$
Again, define $G(\theta)=g_{\eta}(\theta,\eta^*),$ then the corresponding asymptotic equation reads as
$$
G(\theta) = K^{-1}\overline{h}_N(\theta,\eta).
$$
Adapting the idea of the i.i.d. ECF one can easily show that the optimal choice of $K$ is $K=K(\theta)=C(\theta)$ with $C$ being an $M \times M$ matrix with entries
$$
C_{k,l}(\theta) = \E \left[ ~h^*_{k,n}(\theta,\eta^*) h_{l,n}(\theta,\eta^*) \right].
$$
Define the
($\theta$-dependent) asymptotic cost function as
$$W_{E}(\theta,\eta)=\E\left|K^{-1/2}\overline{h}^{(s)}_{n}(\theta,\eta)\right|^2.$$
For each the fixed $\theta$ define $\eta^*(\theta)$ such that
$$W_{E,\eta}(\theta,\eta^*(\theta))=0.$$
Let the Hessian of $W_{E}$ w.r.t. $\eta$ at $\eta = \eta^*(\theta)$ be denoted by
$$R_{E}^{*}(\theta)=W_{E,\eta\eta}(\theta,\eta^*(\theta)).$$
To formulate our result we need some technical conditions. Conditions 1 and 2 have been already presented in Section \ref{sec:disc_levy}.
Let $\rho$ be the joint parameter i.e. $\rho=(\theta,\eta).$ Let $D_{\rho}$ and $D_{\rho}^*$ be compact domains such that \\ $\rho^* \in D_{\rho}^* \subset \text{int } D_{\rho}$ and $D_{\rho} \subset G_{\rho}.$
\begin{cond}\label{cond:3'}
For each $\theta \in D_{\theta}$ the equation $W_{E,\eta}(\theta,\eta)=0$ have a unique solution in $D_{\eta}^*.$
\end{cond}
%\cite{gerencs_arma}.
%{\bf !!! MATE: EXPLAIN NOTATIONS BELOW!}
\begin{lem}\label{lemma:st2_N^-1/2}
Under Conditions \ref{cond:STAB}, \ref{cond:MOMENTS_of_NU} and \ref{cond:3'} we have $\hat{\eta}_N(\theta)-\eta^*(\theta)=O_M(N^{-1/2}).$
\end{lem}
Our next result characterizes the estimation error of the ECF method for the noise parameter $\eta^*.$
\begin{thm}\label{thm:step2_eta_dependent_dif}
Under Conditions \ref{cond:STAB}, \ref{cond:MOMENTS_of_NU} and \ref{cond:3'} we have
$$
\hat{\eta}_N(\theta)-\eta^*(\theta)=-(R_{E}^{*}(\theta))^{-1}V_{E,N \eta
}(\theta,\eta^*(\theta))+O_M(N^{-1}).
$$
\end{thm}
%uj{\bf THE PROOFS SHOULD BE GIVEN OR SKETCHED IN THE APPENDIX}
The proofs of Lemma \ref{lemma:st2_N^-1/2} and Theorem
\ref{thm:step2_eta_dependent_dif} are isomorphic to that of Lemma
\ref{lemmma:PE_N^-1/2} and Theorem \ref{thm:theta_1_diference}.
Using this theorem with $\theta= \hat{\theta}_N,$ the estimation
that we obtained by the PE method, we may conclude the following
result.
\begin{thm}\label{thm:step_2_etaN_dependent_dif_adaptive}
Under Conditions \ref{cond:STAB}, \ref{cond:MOMENTS_of_NU} and \ref{cond:3'} we have
$$
\hat{\eta}_N-\eta^*=-(R_{E}^{*}(\theta^*))^{-1}V_{E,N \eta
}(\theta^*,\eta^*)+O_M(N^{-1}).
$$
\end{thm}
The proof of the last result is obtained by the very same methods
as Lemma \ref{lemmma:PE_N^-1/2} and Theorem
\ref{thm:theta_1_diference}, combining with
$\eta^*(\hat{\theta}_N)-\eta^*(\theta^*)=\eta^*(\hat{\theta}_N)-\eta^*=O_M(N^{-1/2}),$
and that
$$\left\|W^{-1}_{\eta\eta}(\theta^*,\eta^*)-W^{-1}_{\eta\eta}(\hat{\theta}_N,\eta^*(\hat{\theta}_N))\right\|=O_M(N^{-1/2}).$$
%combined with the fact that
%\begin{equation}\label{eq:dif4}
%\left\|W_{\eta\eta}(\theta^*,\eta^*)-W_{\eta\eta}(\hat{\theta}_N,\eta^*)\right\|=O_M(N^{-1/2}),
%\end{equation}
%which is implied by $\hat{\theta}_N-\theta^*=O_M(N^{-1/2}).$ %Equation (\ref{eq:dif4}) and equations (\ref{eq:dif1}), (\ref{eq:dif2}) together imply (\ref{eq:dif3}).
%
\section{Re-estimation of $\theta^*$ by the ECF method for system
identification. The ECF method for identifying the dynamics of a L\'evy systems}
\label{sec:ECF_system}

If we were given the true value of $\eta^*$ the score function would be
$$
h^{opt}_{k,n}(\theta)= \left(e^{iu_k \varepsilon_n (\theta)}-\varphi(u,\eta^*)\right) \varepsilon_{n \theta} (\theta).
$$
Since we are given only an estimation $\hat{\eta}_N$ of $\eta^*$ we will use the score function
\begin{equation}
h_{k,n}(\theta)= \left(e^{iu_k \varepsilon_n (\theta)}-\varphi(u,\hat{\eta}_N)\right) \varepsilon_{n \theta} (\theta).
\end{equation}
%The continuity of $\varphi$ and Lemma \ref{lemma:st2_N^-1/2} together imply that
%\begin{equation}\label{eq:opt_nonopt_diff}
%h^{opt}_{k,n}(\theta)-h_{k,n}(\theta)=O_M(N^{-1/2}).
%\end{equation}
In this section we analyze the identification of $\theta^*$ with an arbitrary given noise characteristic $\eta.$ We prove consistency and we give the precise characterization of the estimation error. As we will see the same results remain valid if we work with $h_{k,n}(\theta)$ instead of $h^{opt}_{k,n}(\theta).$
%One can easily mimic the steps of the proofs and use (\ref{eq:opt_nonopt_diff}) to conclude the same results.
%\section{The ECF method for identifying the dynamics of a L\'evy systems}
%\label{sec:ECF}
The ECF method has been widely used in finance as an alternative to the ML method, assuming i.i.d. returns \cite{CARRASCO-CGMM}, \cite{CARRASCO-EFFEMPIRCHAR}, \cite{YU-EMPIRCHARFUN}. We adapt this technique to the problem of identifying a discrete-time L\'evy system as described in (\ref{eq:disc_levy2}). In this section it is assumed that the characteristic function of the noise, or equivalently $\eta^*$ is known. The problem we address is to identify the system dynamics specified by $\theta^*.$
%%% The study of this problem is motivated by the third step of out identification algorithm.
%%Although there only a good estimation $\hat{\eta}_N$ of $\eta^*$ is available, the result that under certain technical conditions $\hat{\eta}_N-\eta^*=O^{Q/2}_M(N^{-1/2})$ makes the following analysis applicable for last step of the three-stage identification method.
%For discrete time L\'evy systems we define the innovation process by
%
%uj{\bf !!! MINDENT KETSZER ? MOND ???}
%
%
%\begin{equation}
%\label{eq:inverse}
%\varepsilon_n (\theta)=A^{-1}(\theta) \Delta y_n,
%\end{equation}
%with zero initial conditions and $\theta \in D_{\theta},$
%
%uj{\bf !!! OR $\theta \in G_{\theta}.$}
%
%
%Then we have for $n \ge 0$
%\begin{equation}
%\label{eq:statdiff}
%\varepsilon^{(s)}_n(\theta)=\varepsilon_n(\theta)+r_n,
%\end{equation}
%where $r_n=O_M^Q(\alpha^n)$ with some $0<\alpha<1$,
%
%uj{\bf !!! MOVE UP NEXT BIT:
%meaning that for all $1 \leq q <Q$
%$$ \sup_{n} \alpha^{-n} \E^{1/q} \left|r_n\right|^q < \infty.$$}
%explained earlier
Following the philosophy of the ECF method take a fix set $u_i$-s,
$1 \leq i \leq M.$
The first natural candidate for a score function would be
$$ f_{k,n}(\theta,\eta)= e^{iu_k \varepsilon_n
(\theta)}-\varphi(u,\eta),$$
see \cite{ECC_own}.
%uj{\bf !!! WHAT FOLLOWS IS A BIT MESSY. NOW YOU PAY THE PRICE FOR NOT PRESENTING ECF FOR RI.I.D DATA. WHERE IS THE EQUATION}
%$$
%G K^{-1} {\overline h} (\theta) = 0 ???
%$$
%{\bf WHERE IS THE MOTIVATION FOR USING EXACTLY THE IV THAT IS USED ???
%RECHECK MY ZABCZYK SLIDES.}
%
%The derivative of $f$
%with respect to $\theta$ is given by
%$$ e^{iu_k \varepsilon_n (\theta)} iu_k \varepsilon_{n\theta} (\theta),$$
%whose expected value is 0 if $\E[\Delta Z_n]=0.$
%
%uj{\bf !!! WHY IS THIS A PROBLEM ???}
It turns out that the identification method that uses $f$-s as score functions does not give an efficient estimator.
%which may be not a full-rank matrix.
The %$p$-dimensional
%uj{\bf STRESS: $h_{k,n}(\theta) \in {\mathbb R}^p $ !!!}
score functions to be used following the basic idea of the ECF method and the instrumental variable method
are defined as
\begin{equation}
    h_{k,n}(\theta,\eta)= \left(e^{iu_k \varepsilon_n (\theta)}-\varphi(u,\eta)\right) \varepsilon_{n \theta} (\theta)
\end{equation}
with $h_{k,n}(\theta)$ being $p \times 1$ vectors. While $h_{k,n}$
is the function that can be computed in practice, $h^{(s)}_{k,n}$
is easier to handle, because of its stationarity. These are indeed
appropriate score functions, since we obviously have
$$ \E   \left[h^{(s)}_{k,n}(\theta^*,\eta^*)\right]=0$$
%uj{\bf ??? and ??? IS THIS ARGUMENT RELEVANT HERE ???}
%dropped
%
%
% $$ h_{k,n}(\theta^*)=h^{(s)}_{k,n}(\theta^*)+O^Q_M(\alpha^n),$$
%with $0<\alpha<1.$
%To see this write $\E   \left[h^{(s)}_{k,n}(\theta^*)\right]$ as
%\begin{equation}
%\begin{split}
%&\E \left[ \E \left[\left(e^{iu_k \varepsilon^{(s)}_n (\theta^*)}-\varphi(u,\eta^*)\right) \varepsilon^{(s)}_{n \theta}(\theta^*)|\mathscr{F}_{n-1}^{\Delta Z}\right]\right]= \\
%&\E \left[ \varepsilon^{(s)}_{n \theta}(\theta^*)\E \left[\left(e^{iu_k \varepsilon^{(s)}_n (\theta^*)}-\varphi(u,\eta^*)\right)|\mathscr{F}_{n-1}^{\Delta Z}\right]\right]=0.
%\end{split}
%\end{equation}
%{\bf MOVE UP: where $\mathscr{F}_{n-1}^{\Delta Z}=\sigma \left\{ \Delta Z_k: k \leq n-1 \right\},$
%}
%and we used that $\varepsilon^{(s)}_{n \theta}(\theta^*)$ is $\mathscr{F}_{n-1}^{\Delta Z}$ measurable.
%uj{\bf THIS IS A GENERAL REMARK THAT SHOULD APPEAR MUCH EARLIER: While $h_{k,n}$ is the function that can be computed in practice, $h^{(s)}_{k,n}$ is easier to handle, because its stationarity.}
Fix a set of $u$-s: $u_1,\ldots,u_M$ and define
$h_n(\theta)=\left(h^T_{1,n}(\theta),\ldots,h^T_{M,n}(\theta)\right)^T.$
%uj{\bf !!! STRESS AGAIN: the $p$-dimensional sample mean of the $k$-th score as}
Define the $p$-dimensional sample mean of the score vector as:
$$\overline{h}_N(\theta)=\frac{1}{N}\sum_{n=1}^N h_{n}(\theta).$$
Let $K>0$ be a fixed symmetric, positive definite $Mp \times Mp$  $K$
%uj{\bf !!! $Lp \times Lp$ ??? MIT KERES ITT A $K$ ???}
%OK this way
weighting matrix. Define the $Mp$-dimensional $g(\theta,\eta)=\E \left[\overline{h}_N (\theta,\eta)\right].$
%uj{\bf Note that $\theta^*$ is the solution of the over-determined set of non-linear algebraic equations}
Note that $\theta=\theta^*$ is the solution of the over-determined set of non-linear algebraic equations
$$
g(\theta^*,\eta^*)=0. %\quad n=1,\ldots,L
$$
%uj{\bf !!! ADJUST: we seek a least-square solution.}

%uj{\bf !!! HAVE A SECOND THOUGHT ON THIS: Since $g$ is not computable we approximate it by $\overline{h}$ and we define the cost functions as}
Since $g$ is not computable we approximate it by $\overline{h}$ and we seek a least-square solution. Define the cost functions as
\begin{equation}
V_N=V_N(\theta,\eta)=|K^{-1/2}\overline{h}_N(\theta,\eta)|^2,
\end{equation}
%uj{\bf !!! I WOULD NOT REPEAT EACH EQUATION BY PRESENTING ITS STATIONARY APPROXIMATION !}
and obtain $\hat{\hat{\theta}}_N(\eta)$ by solving
\begin{equation}
V_{N \theta}(\theta,\eta)=0.
\label{eq:score2_2}
\end{equation}
The system of equations in (\ref{eq:score2_2}) is no longer over-determined because ${\rm dim } V_{N \theta}=p.$
%uj{\bf !!! WRITE OUT THIS EQUATION !!! GET NEW SCORES AS IN THE I.I.D. CASE !}
This gradient equation can be written as
\begin{equation}
\overline{h}^*_{N \eta}(\theta,\eta) K^{-1} \overline{h}_N(\theta,\eta)=0,
\end{equation}
and this $p$ equations can be considered as a set of new score
functions. The corresponding asymptotic problem can be formulated
as
\begin{equation}
G(\theta) K^{-1} \overline{h}_N(\theta,\eta)=0,
\end{equation}
with $G(\theta)=g_{\eta}(\theta,\eta^*).$
%{\bf !!! WRITE OUT THIS EQUATION !!! GET NEW SCORES AS IN THE I.I.D. CASE !}
%\section{Analysis}
%\label{sec:ECF_proofs}

%uj
%C appears later
The asymptotic cost function is defined by
$$W(\theta,\eta)=\lim_{N \rightarrow \infty} \E[V_{N}(\theta,\eta)]=g^*(\theta,\eta)K^{-1}g(\theta,\eta).$$
Let $\theta^*(\eta)$ denote the $\eta$ dependent the solution of the asymptotic equation
$$W_{\theta}(\theta,\eta)=0.$$
\begin{cond}\label{cond:3''}
For each $\eta \in D_{\eta}$ the equation $W_{\theta}(\theta,\eta)=0$ have a unique solution in $D_{\theta}^*.$
\end{cond}
Note that $\theta^*(\eta)=\theta^*$ for each $\eta$ holds, because $\E\left[ \varepsilon_{n \theta}(\theta^*)\right]=0.$
The Hessian of $W$ at $\theta=\theta^*(\eta)=\theta^*$:
$$R^*(\eta)=g_{\theta}^*(\theta^*,\eta)K^{-1}g_{\theta}(\theta^*,\eta).$$
The following result, which can be proved using the reasoning seen
in \cite{gerencs_arma}, is a martingale representation theorem for
the $\eta$-dependent estimate of $\theta^*.$
%uj{\bf DO WE REALLY WANT TO MONITOR THE DEPENDENCE OF $Q$ ON $p$ AND LATER ON $r$ ALL THE WAY ???}
%no
\begin{thm}\label{thm:step2_theta_dependent_dif}
Under Conditions \ref{cond:STAB}, \ref{cond:MOMENTS_of_NU} and \ref{cond:3''} we have
$$
\hat{\hat{\theta}}_N(\eta)-\theta^*=-(R^{*}(\eta))^{-1}V_{N \theta
}(\theta^*,\eta)+O_M(N^{-1}).
$$
\end{thm}

%uj{\bf !!!  IS THIS REMARK ON ITS RIGHT PLACE ??? For the definition of $L$-mixing processes and for other corresponding definitions and theorems see the Appendix.}
%remark can be found much earlier

%uj{\bf !!! DO WE REALLY NEED TO PRESENT A FULL PROOF THAT IS ISOMORPHIC TO AN OLDER PROOF ???}
%sketch of proof:

\textbf{Sketch of the proof:} First, note that since $\Delta
y_n=\sum_{i=0}^r a_i(\theta^*)\Delta Z_i,$ holds, $\Delta y_n$ is
a linear combination of $L$-mixing processes. Using the fact that
a uniformly exponentially stable filter with $L$-mixing input
produces a uniformly $L$-mixing output \cite{gerencs_mixing} we
get that $\Delta y_n$ is an $L$-mixing process. The innovation
process and its derivatives with respect to $\theta$ can be
written as
\begin{align*}
\varepsilon_n(\theta)&=A^{-1}(\theta) \Delta y_n \\
\varepsilon_{n\theta}(\theta)&=A_{\theta}^{-1}(\theta) \Delta y_n \\
\varepsilon_{n\theta\theta}(\theta)&=A_{\theta \theta}^{-1}(\theta) \Delta y_n.
\end{align*}
Again, since $A^{-1}(\theta)$ and its derivative with respect to
$\theta$ are uniformly exponentially stable we conclude that
the processes $\varepsilon_n(\theta),$ \\
$\varepsilon_{n\theta}(\theta)$ and
$\varepsilon_{n\theta\theta}(\theta)$ are $L$-mixing uniformly in
$\theta$.

The next step is to show that for any given $d>0$ the equation $V_{N \theta}(\theta,\eta)=0$ has a unique solution in $D_{\theta}$ and it is in the sphere $S=\left\{\theta:\left|\theta-\theta^*\right|<d\right\}$ with probability at least $1-O(N^{-s})$ for any $s>0,$ see Lemma 2.3. in \cite{gerencs_arma}.

We have
\begin{equation}\label{eq:taylor}
0=V_{N \theta} \left(\hat{\hat{\theta}}_N,\eta\right)=V_{N \theta} \left(\theta^*,\eta\right)+ \overline{V}_{N\theta \theta}(\eta) \left(\hat{\hat{\theta}}_N-\theta^*,\eta\right),
\end{equation}
where $$ \overline{V}_{N\theta \theta}(\eta)= \int_{0}^{1} V_{N\theta \theta}\left(\left(1-\lambda\right)\theta^*+\lambda \hat{\hat{\theta}}_N,\eta\right) d\lambda .$$
One may proceed like in the proof of Theorem 2.1. in \cite{gerencs_arma} to conclude that
\begin{equation}
\left\|\overline{V}^{-1}_{N\theta\theta}(\eta)-W^{-1}_{\theta\theta}(\theta^*,\eta)\right\|=O_M(N^{-1/2}).
\end{equation}
except from an event of probability $O_M(N^{-s})$ for any $s>0.$
Finally,
\begin{align*}
    \hat{\hat{\theta}}_N-\theta^*=-\overline{V}^{-1}_{N\theta\theta}(\eta)V_{N\theta}(\theta^*,\eta)= \\
    -\left(W^{-1}_{\theta\theta}(\theta^*,\eta)+O_M(N^{-1/2})\right) V_{N\theta}(\theta^*,\eta)= \\
    -W^{-1}_{\theta\theta}(\theta^*,\eta) V_{N\theta}(\theta^*,\eta)+O_M(N^{-1})
\end{align*}
holds, again except from an event of probability $O_M(N^{-s})$ for any $s>0,$ hence the last expression reads as
\begin{equation*}
-\left(R^*(\eta)\right)^{-1}V_{N\theta}(\theta^*,\eta)+O_M(N^{-1}).
\end{equation*}

%\end{pf}

By choosing $\eta$ to be equal to $\hat{\eta}_N,$ the estimate of the noise that we obtained at the second step of the procedure, and using that $R^*$ and $V_{N \theta}$ is smooth enough in $\eta$ and that $\hat{\eta}_N-\eta^*=O_M(N^{-1/2})$ we obtain the following result.
\begin{thm}\label{thm:step2_theta_dependent_dif_adaptive}
Under Conditions \ref{cond:STAB}, \ref{cond:MOMENTS_of_NU} and \ref{cond:3''} we have
$$
\hat{\hat{\theta}}_N(\hat{\eta}_N)-\theta^*=-(R^{*}(\eta^*))^{-1}V_{N \theta
}(\theta^*,\eta^*)+O_M(N^{-1}).
$$
\end{thm}
\section{Efficiency of the single term ECF method}
\label{sec:efficiency}
In this section we compute the asymptotic covariance of the estimator proposed in Section \ref{sec:ECF_system}. Recall that
$$ R_{P}^*=\E \left[ \varepsilon^{(s)}_{n \theta} (\theta^*)\varepsilon^{(s)T}_{n \theta}(\theta^*) \right], $$
and define the $L \times L$ matrix $C$ with elements
$$
C_{k,l} = \varphi(u_k - u_l,\eta^*) - \varphi(u_k ,\eta^*) \varphi(- u_l,\eta^*).
$$
\begin{thm}
Choosing $K=C \otimes R_P^*,$ the inverse of the asymptotic
covariance matrix of the estimator presented in Section
\ref{sec:ECF_system} is
$$
N \left(\psi^* C^{-1} \psi\right)^{-1}
\left(R_P^*\right)^{-1},
$$
where
%uj{\bf !!! $\frac1N$ ???}
$\psi=\left( iu_1 \varphi(u_1),\ldots,iu_M \varphi(u_M)\right)^T.$
\end{thm}
%uj{\bf !!! I WOULD RATHER TAKE THE INVERSE OF $\left(u^* C^{-1} u\right) \Sigma_P^*.$ }

%----
%\begin{pf}
An essential property of $\hat{\eta}_N$ and $\hat{\hat{\theta}}_N(\hat{\eta}_N)$ that they are asymptotically uncorrelated. This can be seen using direct calculation using the fact that $\E\left[ \varepsilon_{n\theta}(\theta^*)\right]=0.$
Using this observation and Theorem \ref{thm:step2_theta_dependent_dif_adaptive} we get that the
%
%uj{\bf !!! ADD ERROR TERMS !}
%
covariance matrix of the estimator is
\begin{equation*}%\nolabel
\begin{split}
&{\rm Cov}\left((\hat{\hat{\theta}}_N-\theta^*)(\hat{\hat{\theta}}_N-\theta^*)^T\right)= \\
&(R^*)^{-1}\E \left[ V^*_{N \theta}(\theta^*)V_{N \theta}(\theta^*)\right](R^*)^{-1}+O_M(N^{-2}).
\end{split}
\end{equation*}
To calculate the above expected value we first approximate $\overline{h}^*_{\theta}(\theta^*)$ with $g_\theta(\theta^*)$:
%uj{\bf !!! ITT AZERT MEGALLNEK EGY PILLANATRA, HA MAR KORABBAN EZ A LEPES PROBLEMSA VOLT, ES ORULNEK ANNAK, HOGY LEHET JOL APPROXIMALNI.}
%Calculating the expected value yields
\begin{equation}%\nolabel
\begin{split}
\E \left[ \overline{h}^*_{\theta}(\theta^*) K^{-1} \overline{h}(\theta^*) \overline{h}^*(\theta^*) K^{-1} \overline{h}_{\theta}(\theta^*)\right]= \\
\E \left[ \left(g^*_\theta(\theta^*)+r^*_1\right)K^{-1} \overline{h}(\theta^*) \overline{h}^*(\theta^*) K^{-1} \left(g_\theta(\theta^*)+r_1\right)\right], \label{eq:cov1}
\end{split}
\end{equation}
with $r_1=O_M^Q(N^{-1/2}).$
Now we calculate the covariance $\E\left[ \overline{h}(\theta^*) \overline{h}^*(\theta^*) \right]$
%uj{\bf !!! EXPLAIN AND STRESS WHAT IS COMING. IT IS GOOD LUCK AND MAGIC !}
$$\E\left[ \overline{h}(\theta^*) \overline{h}^*(\theta^*) \right]= \frac{1}{N} C \otimes R_P^*+r_2=\frac{1}{N}K+r_2,$$ where $r_2=O_M^Q(N^{-1}).$
Using this (\ref{eq:cov1}) reads as
\begin{equation*}
\begin{split}
&\frac{1}{N} g^*_\theta(\theta^*)K^{-1}g_\theta(\theta^*)+\frac{1}{N} \E \left[ r^*_1  K^{-1} r_1\right]+ \\
&\frac{1}{N} \E \left[ g^*_\theta(\theta^*) K^{-1} r_1\right]+
\frac{1}{N} \E \left[ r^*_1  K^{-1} g_\theta(\theta^*)\right]= \\
&\frac{1}{N} g^*_\theta(\theta^*)K^{-1}g_\theta(\theta^*)+r_3,
\end{split}
\end{equation*}
with $r_3=O_M^{Q}(N^{-3/2}),$ because $g^*_\theta(\theta^*)$ is
bounded. Hence, considering that $g^*_\theta(\theta^*)=\psi^*
\otimes R_P^*,$ and using the mixed-product property and the
inverse of a Kronecker product, reading as $(A \otimes B)(C
\otimes D)=AC \otimes BD$ and $(A \otimes B)^{-1}=A^{-1}\otimes
B^{-1},$ the covariance can be calculated as follows:
%uj{\bf !!! EXPLAIN THE SOURCE OF THE EXPRESSION $\left((\psi^* \otimes R_P^*)\right)$ ! ADD A WORD ON TENSOR PRODUCTS, AND ON THE RULE WE APPLY !}
\begin{equation*}%\nolabel
\begin{split}
&{\rm Cov}\left((\hat{\hat{\theta}}_N-\theta^*)(\hat{\hat{\theta}}_N-\theta^*)^T\right)= \\
&\frac{1}{N}\left((\psi^* \otimes R_P^*)(C \otimes R_P^*)^{-1}(\psi \otimes R_P^*)\right)^{-1}+r_3= \\
&\frac{1}{N}\left((\psi^* C^{-1} \psi) \otimes R_P^* \right)^{-1}+r_3= \\
&\frac{1}{N}(\psi^* C^{-1} \psi)^{-1} (R_P^*)^{-1}+r_3,
\end{split}
\end{equation*}
which concludes the proof.
%----
%\end{pf}

\subsection{Efficiency of the estimation procedure}

Now we are ready to demonstrate that the proposed estimation
method is asymptotically efficient. Use the full continuum of
$u$-s and define $K=C$ as an operator like in
\cite{CARRASCO-EFFEMPIRCHAR}
%uj{\bf !!! WRITE OUT IT HERE !}
\begin{equation}
(Cf)(s)=\int c(s,t)f(t) \pi(t) dt,
\end{equation}
with $\pi$ being a probability measure on $\mathbb{R},$ and
$$
c(s,t)=\E\left[ h^*_{s,n}(\theta^*,\eta^*)h_{t,n}(\theta^*,\eta^*)\right],
$$
where the full continuum of $u$-s is defined via $u_s=s$ for all $s \in \mathbb{R}.$

Since $\psi=\left( i\psi_1 \varphi(u_1),\ldots,i\psi_M \varphi(u_M)\right)^T$ if $M$ moment conditions is used, the asymptotic covariance matrix of the estimator with full continuum of $u$-s would be
\begin{equation}\label{eq:fullcont_cov}
\left(|| iu\varphi(u,\eta^*)||_C^2 \right)^{-1} (R_P^*)^{-1}.
\end{equation}
%uj{\bf !!! THIS EXPRESSION DOES NOT QUITE MATCH WHAT WE HAD IN THE LEMMA (WHICH I WOULD RATHER PRESENT AS A THEOREM !)}
%uj{\bf THE FOLLOWING REMARK SHOULD BE REPHRASED.
Note that in the above formula $|| iu\varphi(u,\eta^*)||_C^2$ depends only on the noise characteristics and $R_P^*$ depends on the derivative of the innovation process, hence on the parameters of the linear system.
%Intuitively in the above formula $(R_P^*)^{-1}$ represents the effect of the underlying system and $u\varphi(u,\eta^*)$ represents the noise characteristics.
According to (\ref{eq:mu}) asymptotic efficiency is reached if $$\left(|| iu\varphi(u,\eta^*)||_C^2 \right)^{-1}=\mu,$$ with
$$
\mu=%\lim_{n \rightarrow \infty}
\E \left[\left(\frac{f'(\Delta Z_n,\eta^*)} {f(\Delta
Z_n,\eta^*)}\right)^2\right],$$ where $\mu$
%uj{\bf !!! WAS SHOWN TO BE EQUAL TO ...}
was shown to be equal to the Fisher of the location parameter.

%{\bf !!! WHAT FOLLOWS IS A BIT TOO FAST, AND TOO SUPERFICIAL. THIS MAY ACTUALLY BE STATED IN THE SECTION FOR ECF FOR I.I.D. DATA.}
%agreed to leave it here

According to \cite{CARRASCO-EFFEMPIRCHAR} for i.i.d. samples the
ECF method with continuum $u$-s gives an asymptotically efficient
estimate of an unknown parameter $\lambda^*$ with asymptotic
covariance
$$(||\varphi_{\lambda}(u,\lambda^*)||_C^2)^{-1}.$$

%{\bf !!! WHAT FOLLOWS IS A BIT TOO FAST, AS WELL, AND TOO LONG. MAY ACTUALLY BE STATED IN THE SECTION FOR ECF FOR I.I.D. DATA.}
%agreed to leave it here
Now we show that $\left(|| u\varphi(u,\eta^*)||_C^2 \right)^{-1}$
can be obtained as the asymptotic covariance of an efficient
i.i.d. ECF method, thus the efficiency of our identification
method for L\'evy system follows. Consider the following
identification problem: given a sequence of i.i.d. samples with
distribution $X+\lambda^*$, where $\lambda^*$ is a location
parameter to be estimated, and $X$ is a random variable with known
characteristic function. Let $\varphi_{X+\lambda}$ denote the c.f.
of $X+\lambda,$ then
\begin{equation}
\begin{split}
\frac{\partial}{\partial \lambda}\varphi_{X+\lambda}(u,\eta) = \frac{\partial}{\partial \lambda}\E\left[e^{iu(X+\lambda)}\right]= \\
\frac{\partial}{\partial \lambda}\left( e^{iu\lambda}\E\left[e^{iuX}\right] \right)=iu\varphi_{X+\lambda}(u,\eta),
\end{split}
\end{equation}
%uj{\bf !!! THUS ??? THE EFFICIENT ???}
thus the ECF method for i.i.d samples that estimates $\lambda^*$
gives an asymptotic covariance
$$
\left(|| iu\varphi(u,\eta^*)||_C^2 \right)^{-1},
$$
hence the efficiency follows.
%uj{\bf !!! THIS IS MUCH TOO COLLOQUIAL !}
\section{Discussion}
\label{sec:discussion}

%uj{\bf  !!! THIS IS A CHALLENGE. IF WE ASSUME LESS, THEN WE HAVE NO
%GENERAL THEOREM ! IF WE ASSUME ALMOST NOTHING, THEN THE EFFECT OF
%THE ESTIMATION ERRORS WILL BE A NIGHTMARE ! discussionbe megjegyzes}

%uj{\bf !!! ADD A WORD ON TERDIK SOMEWHERE !}
%done

%uj{\bf IT WOULD BE SILLY TO INDICATE THAT WE HAVE NOT BEEN PRECISE.}

%uj{\bf IT WOULD BE SILLY TO INDICATE THAT WE HAVE NOT BEEN PRECISE.}

We briefly sketch the identification method that uses blocks of dependent data.
Let us consider the parametric family of time series
$$
\Delta y_n(\theta,\eta)=A(\theta)\Delta Z_n(\eta),
$$
with $ - \infty < n < + \infty.$ Note that for $(\theta,\eta) =
(\theta^*,\eta^*)$ we recover our observed data, at least in a
statistical sense.

The method proposed in the literature
%
%\medskip
%\noindent {\bf !!! DO PROVIDE SOME REFERENCES !}
%\medskip
%
is based on the observation that, as an alternative to the joint
probability density function, we can compute the joint
characteristic function of blocks of unprocessed data, i.e. for
blocks of the time series $(y_n).$ Indeed, fix a block length, say
$r,$ and define the $r$-dimensional blocks
$$\Delta Y^r_n(\theta,\eta)=(\Delta y_{n-1}(\theta,\eta),\ldots,\Delta y_{n-r}(\theta,\eta)).$$
%\medskip
%\noindent {\bf !!! WHY FORWARD BLOCKS ??? A MORE STANDARD OPTION
%WOULD BE
%$$\Delta Y^r_n(\theta,\eta)=(\Delta y_{n-1}(\theta,\eta),\ldots,\Delta y_{n-r}(\theta,\eta)).$$
%THEN WE WILL BE INLINE WITH WHAT COMES LATER: THE CONVOLUTION OF U
%AND DELTA Y: THE INDEXING OF DELTA Y AND U WILL BE IN OPPOSITE
%DIRECTIONS.}
%\medskip
Then the characteristic function of $\Delta Y^r_n (\theta,\eta)$,
with $u=(u_1,...,u_r)^T$ being an arbitrary vector in $\mathbb
R^r,$
%\medskip
%\noindent {\bf !!! WHY CAP U ??? IT IS NOT A R.V. !}
%\medskip
is given by
%$$
%\varphi_n(u,\theta,\eta)=\E \left[e^{iu^T \Delta
%Y^r_n(\theta,\eta)}\right] = \E \left[e^{i \sum_{j=1}^r u_j \Delta
%Y_{n-j}(\theta,\eta)}\right].
%$$
%Now this can be explicitly computed, at least in theory, with
%$h_l(\theta),~l=0,1,... $ denoting the impulse responses of the
%system $A(\theta),$ as
%\begin{equation}
%\begin{split}
%%%%\varphi_n(U,\theta,\eta)=\E \left[\exp\{iU^T \Delta Y^r_n(\theta,\eta)\}\right]= \\
%\E \left[\exp \left\{ i\sum_{j=1}^r u_j \sum_{l=0}^{\infty} h_l(\theta) \Delta Z_{n-j-l}(\eta)\right\}\right]= \\
%\E \left[\exp \left\{ i \sum_{k=1}^{\infty}\Delta Z_{n-k} (\eta) \sum_{l \ge 0, j +l=k} u_j h_l(\theta) \right\}\right]. \\
%\end{split}
%\end{equation}
%Writing
%$$
%v_k(\theta) = \sum_{j=1}^{r} h_{k-j}(\theta) u_{j},
%$$
%with $h_l(\theta) =0$ for $l < 0,$ and denoting the characteristic
%function of $\Delta Z_n(\eta)$ (for any $n$) by $\varphi_{\Delta
%Z(\eta)},$ we get
\begin{equation}
\begin{split}
%%%\varphi_n(U,\theta,\eta)=\E \left[\exp\{iU^T \Delta Y^r_n(\theta,\eta)\}\right]= \\
\varphi_n(u,\theta,\eta) = \prod_{k=1}^{\infty} \varphi_{\Delta
Z(\eta)}(v_k(\theta)),
\end{split}
\end{equation}
with $v_k(\theta) = \sum_{j=1}^{r} h_{k-j}(\theta) u_{j}$ and $h_l(\theta),~l=0,1,... $ denoting the impulse responses of the
system $A(\theta).$
%
%{\bf MATE: A $\Delta y_n=0$ for $n<0$ nem biztos, hogy praktikus feltetel, ui, nagy $n$-re igy is ugy is rengeteg tag lesz. A $\Delta Z_{j-l}$ helyett $\Delta Z_{n+j-l-l}$ irando. Szoval en inkabb ezt irnam, sot a $v_j$=ket is expliciten kiirnam:}
%
%\begin{equation}
%\begin{split}
%\varphi_n(U,\theta,\eta)=\E \left[\exp\{iU^T \Delta Y^r_n(\theta,\eta)\}\right]= \\
%\E \left[\exp \left\{ i\sum_{j=1}^r U_j \sum_{l=0}^{\infty} h_l(\theta) \Delta Z_{n+j-l-l}(\eta)\right\}\right]= \\
%\prod_{j=0}^{\infty} \varphi_{\Delta Z(\eta)}(v_j(\theta)),
%\end{split}
%\end{equation}
%
%
%
%Note that for the sequence $v(\theta) = (v_k(\theta))$ we can
%write $v(\theta) = h(\theta) * u$, with $*$ denoting convolution.

Now the ECF method would be defined by fitting this theoretical
characteristic function to the empirical characteristic function,
obtained as the arithmetic mean of the individual scores
$$
h_n(u,\theta,\eta)=e^{iu^T\Delta Y^r_n}-\varphi_n(u,\theta,\eta).
$$
Without going into further details we point out that the weakness
of this approach is that the characteristic function
$\varphi_n(u,\theta,\eta)$ is given in terms of an infinite
product, and hence it is not clear how to use it in actual
computations. Moreover, it is pointed out in the literature that
the above ECF method for dynamic models may be less efficient than the ML method, see \cite{CARRASCO-EFFCGMM}.

%uj{\bf THIS WILL BE THE SUBJECT OF A FORTHCOMING PAPER. }
Furthermore, an interesting problem is to implement and analyze a recursive estimation method for the dynamics and noise characteristics of a L\'evy system, this will be the subject of a forthcoming paper.

%{\bf BTW. IF YOU DO NOT GLUE THE TWO ESTIMATORS OF $\theta$ INTO A SINGLE ENTITY THEN YOU WILL HAVE NO STABILITY PROBLEM WITH THE ODE.}
%checked

%{\bf WE MAY MOVE THE SECTION ON DIRECT ECF TO HERE - IF WE HAVE SPACE. }
%moved here, I made it shorter

%uj{\bf "PROMISE LESS, DELIVER MORE."}

Our aim is to give an identification method for Wiener-Hammerstein models using the basic ideas of the ECF method.

\appendix

\section{$L$-mixing processes}

Let $\theta$ be a $d$-dimensional parameter vector.
\begin{defn}
We say that $x_n(\theta)$ is $M$-bounded if for all $q \geq 1$,
$$ M_q(x)=\sup_{n>0, \theta \in D} \E^{1/q}\left|x_n(\theta)\right|^q < \infty $$
\end{defn}
Define $\mathscr{F}_n=\sigma\left\{e_i: i \leq n \right\}$ and $\mathscr{F}^+_n=\sigma\left\{e_i: i > n \right\}$ where $e_i$-s are i.i.d. random variables.

\begin{defn}
We say that a stochastic process $\left(x_n(\theta)\right)$ is $L$-mixing with respect to $\left(\mathscr{F}_n,\mathscr{F}^+_n\right)$ uniformly in $\theta$ if it is $\mathscr{F}_n$ progressively measurable, M-bounded with any positive $r$ and
\begin{equation*}
%\begin{split}
\gamma_q(r,x)=\sup_{n \geq r, \theta \in D} \E^{1/q} \left|x_n(\theta)-\E\left[x_n(\theta)|\mathscr{F}^+_{n-r}\right]\right|^q,
%\end{split}
\end{equation*}
we have for any $q \geq 1,$
$$\Gamma_q(x)=\sum_{r=1}^{\infty} \gamma_q(r,x) < \infty.$$
\end{defn}

%\begin{thm}\label{thm:ineq}
%Let $(u_n), n \geq 0$ be an $L$-mixing process with $\E u_n=0$ for all $n,$ and let $(f_n)$ be a deterministic sequence. Then we have for all $m \geq 1,$
%\begin{equation*}
%\begin{split}
%&\E^{1/(2m)}\left|\sum_{n=1}^N f_n u_n\right|^{2m} \leq \\
%&C_m \left(\sum_{n=1}^N f^2_n\right)^{1/2} M^{1/2}_{2m}(u) \Gamma^{1/2}_{2m}(u)
%\end{split}
%\end{equation*}
%where $C_m=2(2m-1)^{1/2}.$
%\end{thm}

Define $$\Delta x/\Delta^{\alpha} \theta=\left|x_n(\theta+h)-x_n(\theta)\right|/\left|h\right|^{\alpha}$$
for $n\geq0, \theta\neq\theta+h \in D$ with $0<\alpha \leq 1.$

\begin{defn}
We say that $x_n(\theta)$ is $M$-H\"{o}lder continuous in $\theta$ with exponent $\alpha$ if the process $\Delta x/\Delta^{\alpha} \theta$ is $M$-bounded.
\end{defn}

%Now let us suppose that $(x_n(\theta))$ is measurable, separable, $M$-bounded and $M$-H\"{o}lder in $\theta$
%with exponent $\alpha$ for $\theta \in D.$ The realizations of $(x_n(\theta))$ are continuous in $\theta$ almost
%surely hence $$x^*_n=\max_{\theta \in D_0} \left|x_n(\theta)\right|$$ is well defined for almost all $\omega,$
%where $D_0 \subset \text{int } D$ is a compact domain. Since the realizations of $(x_n(\theta))$ are continuous,
%$x^*_n$ is measurable with respect to $\mathscr{F}.$
%
%\begin{thm}\label{thm:moments}
%Assume that $(x_n(\theta))$ is measurable, separable, $M$-bounded and $M$-H\"{o}lder in $\theta$ with exponent
%$\alpha$ for $\theta \in D.$ Then we have for all positive $q$ and $p/\alpha <s,$
%$$M_q(x^*)\leq C \left(M_{qs}(x)+M_{qs}(\Delta x/\Delta^{\alpha} \theta)\right)$$
%where $C$ depends only on $p, q, s, \alpha$ and $D_0, D.$
%\end{thm}
%Choosing $f_n=1$ and $\alpha=1$ and using Theorem \ref{thm:ineq} and \ref{thm:moments} we obtain

\begin{thm}\label{thm:sup}
Let $(u_n(\theta))$ be an $L$-mixing uniformly in $\theta \in D$ such that $\E u_n(\theta)=0$ for all $n \geq 0, \theta \in D,$ and assume that $\Delta u /\Delta \theta$ is also $L$-mixing uniformly in $\theta, \theta+h \in D.$ Then
\begin{equation}
\sup_{\theta \in D_0} \left|\frac{1}{N} \sum_{n=1}^N u_n(\theta)\right|=O_M(N^{-1/2})
\end{equation}
\end{thm}

\begin{thm}\label{thm:implicit}
Let $D_0$ and $D$ be as above and let \\ $W_{\theta}(\theta), \delta W_{\theta}(\theta), \ \theta \in D \subset \mathbb{R}^p$ be $\mathbb{R}^p$-valued continuously differentiable functions, let for some $\theta^* \in D_0, W_{\theta}(\theta^*)=0,$ and let $W_{\theta \theta}(\theta^*)$ be nonsingular. Then for any $d>0$ there exists positive numbers $d',d''$ such that
\begin{equation}
\left|\delta W_{\theta}(\theta)\right|<d' \text{ and } \left\|\delta W_{\theta \theta}(\theta)\right\|<d''
\end{equation}
for all $\theta \in D_0$ implies that the equation $W_{\theta}(\theta)+\delta W_{\theta}(\theta)=0$ has exactly one solution in a neighborhood of radius $d$ of $\theta^*.$
\end{thm}

\bibliographystyle{plain}        % Include this if you use bibtex
\bibliography{autosam}           % and a bib file to produce the

\begin{thebibliography}{1}
%\bibitem{Black-Scholes}
%F.~Black and M.~Scholes (1973).
%\newblock The pricing of options and corporate liabilities,
%\newblock {\it Journal of Political Economy} 81: 345-637.

%\bibitem{Bollerslev}
%T.~Bollerslev (1986).
%\newblock {Generalized autoregressive
%conditional heteroscedasticity\/},
%\newblock {\it Journal of Econometrics} 31: 307-327.
\bibitem{Bacehelier}
{\sc Bachelier L. (1900).}
{\em Th\'eorie de la sp\'eculation.} Annales Scientifiques de l'\'ecole Normale Sup\'erieure, 3 (17), pp.~21 �86.

\bibitem{black&scholes}
{\sc Black F., Scholes M. (1973).}
{\em The pricing of options and corporate liabilities.} Journal of Political Economy, 81, pp.~345-637.

\bibitem{jacod}
{\sc Jacod J., Shiryaev A.N. (2002).}
{\em Limit theorems for stochastic processes (2. ed.).} Springer

\bibitem{Sato}
{\sc Sato K. (1999).}
{\em L\'evy Processes and Infinitely Divisible Distributions.} Cambridge University Press

\bibitem{terdik}
{\sc Terdik Gy., Woyczynski W.A. (2006).}
{\em Rosinski measures for tempered stable and related Ornstein-Uhlenbeck processes.} Probability and Mathematical Statistics, 26 (2), pp.~213--243.

\bibitem{soderstrom}
{\sc Astrom K., Soderstrom T. (1974).} {\em Uniqueness of the maximum likelihood estimates of the parameters of an ARMA model.} Automatic Control, IEEE Transactions on 19.6, pp.~769-773.

\bibitem{CGMY}
{\sc Carr P., Geman H. , Madan D., Yor M. (2000).} {\em The fine structure of asset returns: an empirical investigation.}
Journal of Business, 75 (2), pp.~305-332.

\bibitem {CARRASCO-EFFCGMM}
{\sc M.~Carrasco,  M.~Chernov,  J.-P. Florens, E.  Ghysels}, {\em Efficient   estimation of general dynamic models with a continuum of moment
  conditions}, in
Journal of Econometrics, 140~(2) (2007), pp.~529--573.

\bibitem{CARRASCO-CGMM}
{\sc Carrasco M., Florens J.-P. (2000).} {\em Generalization of
{GMM} to a continuum of moment conditions.}
Econometric Theory, 16~(06), pp.~797--834.

\bibitem{CARRASCO-EFFEMPIRCHAR}
{\sc Carrasco M., Florens J.-P. (2002).} {\em Efficient {GMM} estimation
using the empirical characteristic function.}
Idei working papers, 140


%\bibitem{KOTCHONI-EFFCHAR}
%{\sc M. Carrasco,  R. Kotchoni}, {\em Efficient estimation using the
%characteristic
%  function}, in
 % ......

\bibitem{LSS_ML_GL+MGY+RZ}%
{\sc Gerencser L., Michaletzky Gy., Reppa Z. (2002).}
{\em A two-step maxumum-likelihood identification
of non-gaussian systems.} Proceedings of the 15th IFAC World Congress


%\bibitem{CGMY}
%{\sc P. Carr, H. Geman, D.B. Madan, M. Yor}, {\em The Fine Structure of Asset Returns: An Empirical Investigation}, in
%The Journal of Business, 75  (2000), pp.~305--322.
%{\bf REF:  P. Carr, H. Geman, D.B. Madan, M. Yor: The fine
%structure of asset returns: an empirical investigation.  The
%Journal of Business, {\bf 75} 2002. pp. 305-332.}

\bibitem{raible}
{\sc Raible S. (2000).} {\em L\'evy Processes in Finance: Theory, Numerics, and Empirical Facts.} PhD dissertation, \\ http://www.freidok.uni-freiburg.de/volltexte/51/pdf/511.pdf, Accessed on 19 November 2012

\bibitem{CONT-TANKOV-FinModJunp}
{\sc Cont R., Tankov P. (2006).} {\em  Financial Modelling with Jump Processes.}
Journal of the American Statistical Association, 101, pp.~1315-1316.

%\bibitem{DESOER-SLOW}
%  C.~A. Desoer, ``Slowly time varying system $\dot x = a (t) x$.'' \emph{{IEEE}
%    Trans. Autom. Control}, vol.~14, no.~6, pp. 780---781, 1969.

%\bibitem{DEVROYE-RNDGEN}
% L. Devroye (1986).
%\newblock Non-uniform random variate generation.


%\bibitem{DEWOLF+WIBERG-ODE}
%D.~DeWolf, D.~Wiberg (1993).
%\newblock An ordinary differential equation
%technique for
%  continuous-time parameter estimation,
%\newblock {\it {IEEE} Trans. Autom. Control},
%  vol.~38:  514--528.


%\bibitem{DUAN-1995}
%J.-C.~Duan (1995).
%\newblock The GARCH option pricing model,
%\newblock {\it Mathematical Finance} 5(1): 13-32.

\bibitem{FEUER_JASA}
{\sc Feuerverger A., McDunnogh P. (1981).} {\em On some Fourier methods for inference.} Journal of the American Statistical Association, 76~(374), pp. 379--387.


\bibitem{FEUER+MC}
{\sc Feuerverger A., McDunnogh P. (1981).} {\em On the efficiency of empirical characteristic function procedures.}
J.R. Stat. Soc. B, 43~(1), pp. 20--47.

%\bibitem{G+GY+M}
%L.~Gerencs\'er, I.~Gy\"ongy,, G.~Michaletzky
%(1984).
%\newblock Continuous-time recursive
%  maximum-likelihood method. {A} new approach to {L}jung's scheme,
%\newblock \emph{A
%  Bridge Between Control Science, Technology. Proc. of the 9th Triennal
%  World Congress of {IFAC}, Budapest}, L.~Ljung, K.~{\AA}str\"om, Eds.
%  \newblock Pergamon Press, Oxford, pp. 75--77.


%\bibitem{GL+PROKAJ-RML-ECC09}
%L.~Gerencs\'er, V.~Prokaj (2009).
%\newblock  Recursive identification of
%continuous-time
%  linear stochastic systems - convergence w.p.1, in {$L_q$}, in
%\emph{European Control Conference 2009 ECC'09 23-26 August 2009,
%Budapest,
%  Hungary}, pp. 1209--1214, cD-ROM MoC1.5.

%\bibitem{HYBRID}
%L.~Gerencs\'er, V.~Prokaj (2010).
%\newblock  Stability of a class of hybrid linear stochastic systems,
%\newblock  \emph{{IEEE} Trans. Autom. Control} 55 (5): 1233-1238.

%\bibitem{MENN-RACHEV-2005}
%C.~Menn, S.T.~Rachev (2005).
%\newblock A GARCH option pricing model with $\alpha$-stable innovations,
%\newblock {\it European Journal of Operational Research} 163: 201-209.

\bibitem{CARR-VG}
{\sc Madan B., Carr P., Chang C. (1998).} {\em The Variance Gamma Process and Option Pricing.}
European Finance Review, 2, pp.~79--105.

\bibitem{MANDELBROT}
{\sc Mandelbrot B. (1963).} {\em The Variation of Certain Speculative Prices.} Journal of Business, 35

\bibitem{duncan}
{\sc Duncan T. E. (2004).} {\em Some processes associated with a fractional Brownian motion.} Mathematics of Finance, (eds. G. Yin, Q. Zhang) Contemp. Math. 351, pp.~93--102.


%\bibitem{Kim-Rachev-Bianchi-Fabozzi}
%Y.S.~Kim, S.T.~Rachev, M.L.~Bianchi, F.J.~Fabozzi (2008).
%\newblock Financial Market Models with L\'evy Processes, Time-Varying Volatility,
%\newblock {\it Journal of Banking, Finance} 32: 1363-1378.

%\bibitem{Kim-Rachev-Chung-Bianchi}
%Y.S.~Kim, S.T.~Rachev, D.M.~Chung, M.L.~Bianchi (2009).
%\newblock The modified tempered stable distribution, GARCH models, option pricing,
%\newblock {\it Probability, Mathematical Statistics} 29(1): 91-117.

%\bibitem{Kim-Rachev-Bianchi-Fabozzi-2}
%Y.S.~Kim, S.T.~Rachev, M.L.~Bianchi, F.J.~Fabozzi (2010).
%\newblock Tempered stable, tempered infinitely divisible GARCH models,
%\newblock {\it Journal of Banking, Finance, to appear.}


%\bibitem{LEVANONY+ZEITOUNI-RML}
%D.~Levanony, A.~Shwartz,, O.~Zeitouni (1994).
%\newblock Recursive identification in
%  continuous-time stochastic processes,
%\newblock \emph{Stochastic Process. Appl.},
%49(2): 245--275.

%\bibitem{POIROT+TANKOV}
%J.~Poirot,, P.~ Tankov (2006).
%\newblock Monte Carlo Option Pricing for Tempered Stable (CGMY) Processes,
%\newblock \emph{Asia-Pacific Financial Markets},
%13 (4): 327-344.
%DOI: 10.1007/s10690-007-9048-7

%\bibitem{TORMA-DDMCP}
%B. Torma, B. G.-T\'{o}th (2010).
%\newblock An efficient descent direction
%method with cutting planes.
%\newblock {\it Central European Journal of Operations Research},
%  forthcoming 18.
\bibitem{knight}
{\sc Knight J. L., Yu J. (2002).} {\em Empirical characteristic function in time series estimation.} Econometric Theory, 18~(3), pp.~691--721.


\bibitem{YU-EMPIRCHARFUN}
{\sc Yu J. (2004).} {\em Empirical characteristic function estimation and
its   applications.}
Econometric Reviews, 23~(2), pp.~93--123.

\bibitem{novikov}
{\sc Miyahara Y., Novikov A. (2001).} {\em Geometric L\'evy Process Pricing Model.} Research Paper Series 66, Quantitative Finance Research Centre, University of Technology, Sydney

\bibitem{gerencs_arma}
{\sc Gerencs\'er L. (1990).} {\em On the martingale approximation of the estimation error of ARMA parameters.}
System \& Control Letters, 15, pp.~417--423.

\bibitem{calder_davis}
{\sc Calder M., Davis R.A.  (1997).} {\em Inference for linear processes with stable noise.}
A practical guide to heavy tails, Birkhauser Boston Inc., pp.~159--176.

\bibitem{gerencs_mixing}
{\sc Gerencs\'er L. (1989).} {\em On a class of mixing processes.} Stochastics, 26, pp.~165--191.

\bibitem{gerencs_ML}
{\sc Gerencs\'er L., Michaletzky Gy., Reppa Z. (2002).} {\em A two-step maximum-likelihood identification of non-Gaussian systems.}
Proceedings of the 15th IFAC World Congress, 15

\bibitem{non-neg_carma}
{\sc Brockwell P. J., Davis R. A., Yang Y. (2011).}{\em Estimation for Non-Negative Lévy-Driven
CARMA Processes.} Journal of Business \& Economic Statistics, 29, pp.~250--259.

\bibitem{Brockwell-Schlemm}
{\sc Brockwell P. J., Schlemm E. (2012).} {\em Parametric estimation of the driving Lévy process of multivariate CARMA processes from discrete observations.} Journal of Multivariate Analysis

\bibitem{quasi_max}
{\sc Schlemm E., Stelzer R. (2012).} {\em Quasi Maximum Likelihood Estimation for Strongly Mixing State Space Models and Multivariate CARMA Processes.} Electronic Journal of Statistics, 6, pp.~2185--2234.

\bibitem{wind}
{\sc Benth F. E., Benth S. J. (2009).} {\em Dynamic pricing of wind futures.}
Energy Economics, 31, pp.~16--24.

\bibitem{queueing}
{\sc  Govil M. K., Fu M. C. (1999).} {\em Queueing theory in manufacturing: A survey.} Journal of manufacturing systems, 18, 214.

\bibitem{ljung}
{\sc Ljung L. (1998).} {\em System Identification: Theory for the User.} Pearson Education

\bibitem{ECC_own}
{\sc Gerencs\'er L., M\'anfay M. (2013).} {\em Identification of finite dimensional linear stochastic systems driven by L\'evy processes.}
Proceeding of European Control Conference, pp.~2415--2420.

\bibitem{tien}
{\sc Tien D.V. (2011).} {\em Multi-Server Markov Queueing Models: Computational Algorithms and ICT Applications.} Dissertation for the Doctor degree of
the Hungarian Academy of Sciences

\bibitem{kingman}
J.F.C. Kingman: {\it Poisson Processes}, Oxford Science Publications,
Clarendon Press, Oxford, 1995.


\end{thebibliography}
                                 % bibliography (preferred). The
                                 % correct style is generated by
                                 % Elsevier at the time of printing.

%\begin{thebibliography}{99}     % Otherwise use the
                                 % thebibliography environment.
                                 % Insert the full references here.
                                 % See a recent issue of Automatica
                                 % for the style.
%  \bibitem[Heritage, 1992]{Heritage:92}
%     (1992) {\it The American Heritage.
%     Dictionary of the American Language.}
%     Houghton Mifflin Company.
%  \bibitem[Able, 1956]{Abl:56}
%     B.~C.~Able (1956). Nucleic acid content of macroscope.
%     {\it Nature 2}, 7--9.
%  \bibitem[Able {\em et al.}, 1954]{AbTaRu:54}
%     B.~C. Able, R.~A. Tagg, and M.~Rush (1954).
%     Enzyme-catalyzed cellular transanimations.
%     In A.~F.~Round, editor,
%     {\it Advances in Enzymology Vol. 2} (125--247).
%     New York, Academic Press.
%  \bibitem[R.~Keohane, 1958]{Keo:58}
%     R.~Keohane (1958).
%     {\it Power and Interdependence:
%     World Politics in Transition.}
%     Boston, Little, Brown \& Co.
%  \bibitem[Powers, 1985]{Pow:85}
%     T.~Powers (1985).
%     Is there a way out?
%     {\it Harpers, June 1985}, 35--47.

%\end{thebibliography}

%\appendix
%\section{Proofs}    % Each appendix must have a short title.
%\section{Some Latin vocabulary}         % Sections and subsections are supported
                                        % in the appendices.
\end{document}